\documentclass[a4paper]{amsart}

\usepackage[labelfont=bf, textfont=normal,
			justification=justified,
			singlelinecheck=false]{caption} % settings for caption in figures
\usepackage{graphicx} % to include figures, i.e. 'includegraphics' command
\usepackage{bbm} % used to write bold "1" as an indicator function
\usepackage[top=2.9cm, bottom=2.9cm,
			inner=3cm, outer=3cm]{geometry}			% Page margin lengths
\usepackage{parskip} % creates skip between every paragraph

\usepackage{enumerate}

\usepackage{booktabs}

\usepackage{pgfplots}
\pgfplotsset{compat=1.16}

\usepackage{subcaption}
\usepackage{pgf,tikz}
\usepgfplotslibrary{groupplots}
\usetikzlibrary{spy, external}
\usepackage{caption}
\usepackage{todonotes}

\usetikzlibrary{arrows.meta, positioning, shadows, calc}

\usepackage{float}

\newcommand{\triple}{{\vert\kern-0.25ex\vert\kern-0.25ex\vert}}

\usepackage{comment}
\newcommand*\dd{\mathop{}\!\mathrm{d}}
\newcommand{\R}{\mathbb{R}}

\newcommand*{\E}{\mathbb{E}}

\newcommand{\N}{\mathcal{N}}

\newtheorem{theorem}{Theorem}[section]

\theoremstyle{definition}

\usepackage[
    backend=biber,
    style=numeric,
    sorting=nyt,
    giveninits=true,
    language=english,
    autolang=other
]{biblatex}

\DeclareFieldFormat{sentencecase}{\MakeSentenceCase*{#1}}

\renewbibmacro*{title}{%
  \ifboolexpr{
    test {\iffieldundef{title}}
    and
    test {\iffieldundef{subtitle}}
  }
    {}
    {\ifboolexpr{
       test {\ifentrytype{article}}
       or test {\ifentrytype{inbook}}
       or test {\ifentrytype{incollection}}
       or test {\ifentrytype{inproceedings}}
       or test {\ifentrytype{inreference}}
       or test {\ifentrytype{patent}}
       or test {\ifentrytype{thesis}}
       or test {\ifentrytype{unpublished}}
     }
       {\printtext[title]{%
          \printfield[sentencecase]{title}%
          \setunit{\subtitlepunct}%
          \printfield[sentencecase]{subtitle}}}
       {\printtext[title]{%
          \printfield[titlecase]{title}%
          \setunit{\subtitlepunct}%
          \printfield[titlecase]{subtitle}}}%
     \newunit}%
}

\makeatletter
\newcommand{\printaddresseshere}{%
  \@setaddresses
}
\makeatother

\usepackage{hyperref}

\newcommand{\defeq}{\mathrel{:\mkern-0.25mu=}}

\definecolor{crimson2143940}{RGB}{214,39,40}
\definecolor{darkgray176}{RGB}{176,176,176}
\definecolor{darkorange25512714}{RGB}{255,127,14}
\definecolor{darkturquoise23190207}{RGB}{23,190,207}
\definecolor{forestgreen4416044}{RGB}{44,160,44}
\definecolor{goldenrod18818934}{RGB}{188,189,34}
\definecolor{gray127}{RGB}{127,127,127}
\definecolor{lightgray204}{RGB}{204,204,204}
\definecolor{mediumpurple148103189}{RGB}{148,103,189}
\definecolor{orchid227119194}{RGB}{227,119,194}
\definecolor{sienna1408675}{RGB}{140,86,75}
\definecolor{steelblue31119180}{RGB}{31,119,180}

\definecolor{OIBlue}{RGB}{0,114,178}
\definecolor{OISkyBlue}{RGB}{86,180,233}
\definecolor{OIGreen}{RGB}{0,158,115}
\definecolor{OIOrange}{RGB}{230,159,0}
\definecolor{OIVermilion}{RGB}{213,94,0}
\definecolor{OIPurple}{RGB}{204,121,167}

\pgfplotsset{
  compat=1.18,
  every axis plot/.append style={line join=round},
}

\tikzset{with marks/.style={mark repeat=4, mark phase=2, mark size=1.5pt}}

\tikzset{with marks long/.style={mark repeat=40, mark phase=20, mark size=1.5pt}}

\tikzset{with marks every/.style={mark repeat=1, mark phase=1, mark size=1.5pt}}

\tikzset{method LogBSDEF/.style={color=OIBlue, very thick, mark=*, mark options={solid}}}
\tikzset{method BSDEF/.style={color=OISkyBlue, very thick, mark=o, mark options={fill=white,solid}}}

\tikzset{method LogDSF/.style={color=OIPurple, very thick, mark=square*, mark options={solid}}}
\tikzset{method DSF/.style={color={OIPurple!70!black}, very thick, mark=square, mark options={fill=white,solid}}}

\tikzset{method PF 1e6/.style={color={OIVermilion!95!black}, thick, solid,  mark=triangle*, mark options={solid}}}
\tikzset{method PF 1e5/.style={color={OIVermilion!80!black}, thick, dashed, mark=triangle*, mark options={solid}}}
\tikzset{method PF 1e4/.style={color={OIVermilion!60!black}, thick, dotted, mark=triangle*, mark options={solid}}}

\tikzset{method EnKF 1e6/.style={color={OIGreen!90!black}, thick, solid,  mark=diamond*, mark options={solid}}}
\tikzset{method EnKF 1e5/.style={color={OIGreen!70!black}, thick, dashed, mark=diamond*, mark options={solid}}}
\tikzset{method EnKF 1e4/.style={color={OIGreen!55!black}, thick, dotted, mark=diamond*, mark options={solid}}}

\tikzset{method EKF/.style={color={OIOrange!90!black}, thick, dashdotted, mark=star, mark options={solid}}}
\tikzset{method Reference/.style={color=black!60, ultra thick, mark=none}}

\tikzset{method KF density/.style={color={OIGreen!55!black}, thick, dashed}}

\tikzset{method PF 1e4 density/.style={color={OIVermilion!60!black}, thick, dotted}}

\tikzset{method EnKF 1e4/.style={color={OIGreen!55!black}, thick, dotted, mark=diamond*, mark options={solid}}}

\tikzset{method LogDSF density/.style={color=OIPurple, thick}}

\tikzset{method LogBSDEF density/.style={color=OIBlue, very thick, dotted}}
\tikzset{method LogBSDEF-LSTM density/.style={color=OISkyBlue, thick, dashdotted}}

\title[High-dimensional Bayesian filtering through deep density approximation]{
High-dimensional Bayesian filtering through deep density approximation}

\begin{document}

%\begin{bibunit}
\begin{refsection}

\author[K.~B{\aa}gmark]{Kasper B{\aa}gmark}
\address{Kasper B{\aa}gmark\\
Department of Mathematical Sciences\\
Chalmers University of Technology and University of Gothenburg\\
SE--412 96 Gothenburg\\
Sweden}
\email{bagmark@chalmers.se}

\author[F.~Rydin]{Filip Rydin}
\address{Filip Rydin\\
Department of Electrical Engineering\\
Chalmers University of Technology\\
SE--412 96 Gothenburg\\
Sweden}
\email{filipry@chalmers.se}

\keywords{Filtering problem, high-dimensional, Fokker--Planck equation, backward stochastic differential equations, splitting, deep learning}
\subjclass[2020]{60G25, 60G35, 62F15, 62G07, 62M20, 65C30, 65M75, 68T07}

\begin{abstract}
In this work, we systematically benchmark two recently developed deep density methods for nonlinear filtering. We model the filtering density of a discretely observed stochastic differential equation through the associated Fokker--Planck equation, coupled with Bayesian updates at discrete observation times. The two filters: the deep splitting filter and the deep backward stochastic differential equation filter, are both based on Feynman--Kac formulas, Euler--Maruyama discretizations and neural networks. The two methods are extended to logarithmic formulations providing sound, robust, and positivity-preserving density approximations in increasing state dimension. Comparing to the classical bootstrap particle filter and an ensemble Kalman filter, we benchmark the methods on numerous examples. In the low-dimensional examples the particle filters work well, but when we scale up to a partially observed $100$-dimensional Lorenz-96 model, the particle-based methods fail and the logarithmic deep backward stochastic differential equation filter prevails. In terms of computational efficiency, the deep density methods reduce inference time by roughly two to five orders of magnitude relative to the particle-based filters.
\end{abstract}
\maketitle
\section{Introduction}
\label{Section: Introduction}
Filtering problems play a central role in applications, where hidden stochastic variables must be inferred from noisy and partial observations, with examples ranging from molecular dynamics and signal processing to finance and climate modeling \cite{barshalom2001estimation, blackman1999design, cui2005comparison,Duc_Kuroda,galanis2006applications, goodman1997mathematics, johannes2009mcmc,kamino2023optimal}. A mathematically rigorous description is provided by the filtering equations, which characterize the evolution of the conditional distribution of the hidden state given the observations.

Classical solution strategies include Kalman-type methods and particle filters. While successful in low dimensions, these approaches face severe challenges in high-dimensional and nonlinear settings: the curse of dimensionality causes particle methods to degenerate \cite{bickel2008sharp,snyder2011particle,snyder2015performance}, and linearization-based approximations commonly utilized in Kalman filters lose accuracy \cite{ehrendorfer2007review,HAN20081434}. To overcome these limitations, recent research has explored probabilistic numerical methods that reformulate the filtering problem through associated Partial Differential Equations (PDEs) and Backward Stochastic Differential Equations (BSDEs). Pioneering work has shown that deep neural networks can be used to approximate high-dimensional PDEs and BSDEs with remarkable success \cite{Arnulf_PDE, E_2017}. 

Building on this foundation, the deep splitting method and the deep BSDE method were recently extended and analyzed in the filtering context, leading to the Deep Splitting Filter (DSF) \cite{baagmark2024convergent} and the deep BSDE Filter (BSDEF) \cite{baagmark2025nonlinear}. 
In the present paper, we complement the initial numerical evaluations of the DSF and BSDEF by providing a systematic numerical comparison of the these two approaches. We also propose log-scale variants that solve the corresponding log-equations. Recasting the equations to log-scale yields positivity-preserving schemes for the density approximation and improves numerical stability in regimes where the filtering density becomes small. Across a broad suite of linear and nonlinear examples, spanning moderate to high dimensions, we investigate accuracy, robustness, and computational cost, and benchmark against the Ensemble Kalman Filter (EnKF), the Extended Kalman Filter (EKF), and a bootstrap Particle Filter (PF). Our aim is to examine the regimes in which these deep density approximations, particularly their log variants, offer practical advantages over classical alternatives. 
The code for the implementation is publicly available\footnote{\url{https://github.com/bagmark/deep-density-filtering}}.

\subsection{Contributions}
\label{subsec: contributions}
This paper considers two density based neural approaches to nonlinear filtering, the deep splitting filter and the deep BSDE filter, together with their logarithmic counterparts (LogDSF, LogBSDEF). The log-formulations are derived from the log-transformed filtering equations (Theorem~\ref{theorem: log equation}) and yield positivity-preserving density approximations together with an alternative objective function that avoids vanishing densities and gradients.
All methods are implemented within a common simulation framework, combining Euler--Maruyama discretization for prediction, Bayes’ rule for updates, and normalization through importance sampling with unconditional-Gaussian or EKF-based proposals. Their performance is evaluated on a broad set of linear and nonlinear test cases: long-horizon Ornstein--Uhlenbeck processes in dimension $d=10$, a nonlinearly observed Schlögl model, and the stochastic Lorenz-96 system with dimension increasing up to $d=100$. Comparisons are made against classical baselines including the EKF, EnKF, and PF, using up to four different error metrics. This includes a first moment error, a mean absolute error, the Kullback--Leibler divergence, and the negative log-likelihood of the state. 

To summarize, the main contributions of this work are:
\begin{enumerate}
    \item A log-density reformulation of the filtering equations (Theorem~\ref{theorem: log equation}) and the resulting LogDSF and LogBSDEF formulations, yielding positivity-preserving schemes for the density approximation and providing a theoretically grounded and numerically robust framework.
    \item An extensive numerical comparison of DSF, BSDEF, and their logarithmic variants across low-, moderate-, and high-dimensional filtering problems, with baselines provided by the EKF, EnKF, and PF.  
    \item A demonstration that the LogBSDEF formulation remains accurate and stable in high-dimensional, nonlinear, and partially observed regimes, including satisfactory performance on the chaotic Lorenz-96 system with dimension $d=100$.
\end{enumerate}

\subsection{Outline}
\label{subsec: outline}
Section~\ref{Section: Background} formulates the filtering problem and the associated PDE system, and establishes the log-density equation used throughout the paper. Section~\ref{Section: Methods} reformulates the DSF and BSDEF frameworks together with their logarithmic variants, including their common probabilistic backbone through Feynman--Kac formulas. Section~\ref{section: experiment setup} summarizes the experimental setup, including evaluation metrics, with further details provided in the supplementary materials. Section~\ref{Section: Experiments} presents the main results, covering numerical experiments for one-dimensional toy examples, high-dimensional linear models, the Schlögl system, and the stochastic Lorenz-96 system. The section ends with a comparison of computational efficiency of the considered methods. 
Section~\ref{section: conclusion} offer a conclusion and a brief discussion. Appendix~\ref{appendix: Proofs} contains the proof of Theorem~\ref{theorem: log equation}. 
The Supplementary Materials provide all details for reproducibility. It is organized into the following parts: Classical methods, Additional example: Linear spring-mass, Normalization, Model architectures, Training, Evaluation, Additional figures, The chemical Schl\"ogl model. % classical baseline methods, a linear spring-mass example, architectures and normalization, training and evaluation details, additional trajectory/density figures, and further details on the Schl\"ogl model.

\subsection{Notation}
\label{subsec: notation}
Throughout this work, $T>0$, is the time horizon and $K\in\mathbb{N}$ the number of observation times, with the uniform grid $t_k = k\frac{T}{K}$, for $k=1,\dots,K$. The state space is $\R^d$ and the observations are in $\R^{d'}$, where $d,d'\in\mathbb{N}$. For an observation sequence we write $o_{1:k}=(o_1,\dots,o_k)\in\R^{d'\times k}$ and denote the observation space $\mathbb{O}_k\defeq\R^{d'\times k}$. The Euclidean norm is $\|\cdot\|$. For scalar fields $\varphi\colon\R^d\to\R$, $\nabla\varphi$ denotes the gradient; for vector-valued maps $g\colon\R^d\to\R^{d'}$, $\mathrm{D} g$ denotes the Jacobian. Function spaces follow standard conventions: $C^k(\R^d;\R^n)$ consists of all $k$-times continuously differentiable functions from $\R^d$ to $\R^n$, $C^{1,2}([s,t]\times\R^d;\R^n)$ the functions which are once differentiable in time and twice differentiable in space, and $L^\infty(\mathbb{X};\mathbb{Y})$ contains the essentially bounded measurable functions $\mathbb{X}\to\mathbb{Y}$. The Gaussian law with mean $m\in \R^d$ and covariance $Q\in \R^{d\times d}$, is $\N(m,Q)$, where we let $\N(x \mid m,Q)$ denote the Gaussian density evaluated at $x\in\R^d$. 

\section{The filtering problem}
\label{Section: Background}
Let $(\Omega,\mathcal{A},(\mathcal{F}_t)_{0\leq t\leq T},\mathbb{P})$ be a complete filtered probability space with the filtration $\mathcal{F} \defeq (\mathcal{F}_t)_{0\leq t\leq T}$ being generated by the two independent $d$-dimensional Brownian motions $B$ and $W$. Assuming sufficient regularity for the functions $\mu \colon \R^d \to \R^d$ and $\sigma \colon \R^d\to \R^{d\times d}$, we have a unique strong solution $S=(S_t)_{t\in[0,T]}$ to the Stochastic Differential Equation (SDE) \cite{Oksendal}
\begin{equation}
    \label{eq: state}
    S_t
    =
    S_0
    +
    \int_0^t
    \mu
    (S_s)
    \dd s
    +
    \int_0^t
    \sigma
    (S_s)
    \dd B_s
    , \quad
    t\in [0,T],
\end{equation}
where $S_0$ is distributed according to some initial distribution $\pi_0$. We let $h\colon \R^d \to \R^{d'}$ be a measurement function, implicitly defining the stochastic observation process $O=(O_k)_{k=1}^K$ by
\begin{align}
    \label{eq: obs}
    O_k
    &=
    h(S_{t_k})
    +
    V_k
    ,\quad
    k=1,\dots,K
    ,
\end{align}
where $V_k\sim \N(0,R)$ with some $d'\times d'$-covariance matrix $R$. From \eqref{eq: obs} we obtain the likelihood function $L(o,x) = \N(o \mid h(x),R)$, which will be used throughout this paper. The methods presented in Section~\ref{Section: Methods} are, however, applicable to more general likelihood models beyond the Gaussian case. Within this setting, our objective is to infer the latent state $S$ given the sequence of observations $O$. More precisely, for $k=1,\dots,K$, and $o_{1:k}\in\mathbb{O}_k$, we want to find the conditional probability density $p_k$ that satisfies, for a measurable set $B\subset \R^d$,
\begin{align*}
    \mathbb{P}
    (S_{t_k} \in B 
    \mid 
    o_{1:k}
    )
    &=
    \int_B
    p_k
    (x
    \mid 
    o_{1:k}
    )
    \dd x
    .
\end{align*}
The conditional density, commonly referred to as the filtering density, satisfies a system of PDEs. To introduce this system we first define $a= (a_{ij}) \defeq  \sigma\sigma^\top$ with entries $a_{ij}$, the infinitesimal generator $A$, defined by the coefficients of \eqref{eq: state}, and its adjoint $A^*$. For $\varphi \in C^2$, $A$ and $A^*$ are given by
\begin{align*}
    A\varphi 
    = 
    \frac{1}{2}\sum_{i,j=1}^d a_{ij}\,
    \frac{\partial^2 \varphi}{\partial x_i \partial x_j} 
    + \sum_{i=1}^d \mu_i \, 
    \frac{\partial \varphi}{\partial x_i}
    ,\quad 
    A^* \varphi 
    = 
    \frac{1}{2}\sum_{i,j=1}^d 
    \frac{\partial^2}{\partial x_i \partial x_j} 
    (a_{ij}\varphi) 
    - \sum_{i=1}^d 
    \frac{\partial}{\partial x_i} 
    (\mu_i \varphi).
\end{align*}
Expanding the derivatives in $A^*$, we see that $(A^*\varphi)(x) = A\varphi(x) + f(x,\varphi(x),\nabla \varphi(x))$, where the function $f\colon \R^d\times\R\times \R^d \to \R$ is defined by 
\begin{align*}
\begin{split}
    f(x,u,v) 
    &= 
    \sum_{i,j=1}^d 
    \frac{
    \partial a_{ij}(x)
    }
    {
    \partial x_i 
    }
    \, v_j
    + 
    \frac{1}{2}
    \sum_{i,j=1}^d 
    \frac{
    \partial^2 
    a_{ij}(x)
    }
    {
    \partial x_i 
    \partial x_j
    }
    \,  u 
    - 
    \sum_{i=1}^d 
    \frac{
    \partial \mu_i(x)
    }
    {
    \partial x_i
    }
    \, u
    - 
    2
    \sum_{i=1}^d
    \mu_i(x) 
    \, v_i
    .
\end{split}
\end{align*}
Assuming sufficient regularity on $\mu$ and $\sigma$, we denote by $p=(p_{k})_{k=0}^K$ the conditional density that solves the unnormalized filtering equations. These equations are initialized at $p_0(0) = \pi_0$, and are recursively given, for $k=0,\dots,K-1$, $x\in\R^d$, $o_{1:k}\in \mathbb{O}_k$, and $t\in[t_k,t_{k+1}]$, by
\begin{align}
    \begin{split}
    \label{eq: FP-prediction}
    &\hspace{1em}p_{k}(t,x,o_{1:k}) 
    = 
    p_{k}(t_{k},x,o_{1:k})
    +
    \int_{t_{k}}^{t} 
    \Big(
    A p_{k}(s,x,o_{1:k}) 
    +
    f
    (
    x,
    p_{k}(s,x,o_{1:k}),
    \nabla 
    p_{k}(s,x,o_{1:k})
    )
    \Big)
    \dd s
    ,
    \end{split}
    \\
    \label{eq: FP-update}
    &\hspace{1em}p_{k+1}
    (t_{k+1},x,o_{1:k+1})
    =
    p_{k}
    (t_{k+1},x,o_{1:k})
    L(o_{k+1},x)
    .
\end{align}
The first equation \eqref{eq: FP-prediction} governs the deterministic evolution of the distribution of $S$ between two observations, while \eqref{eq: FP-update} incorporates the observations through Bayes' formula. The solution $p=(p_k)_{k=0}^K$ consists of $K+1$ piecewise continuous conditional densities satisfying, $p_k(o_{1:k})\in C^{1,2}([t_k,t_{k+1}]\times\R^d)$ for all $k=0,\dots,K-1$ and $o_{1:k}\in \mathbb{O}_k$. This is shown under a parabolic Hörmander condition in \cite{baagmark2024convergent}. 

Considering that $h$ and $R$ are assumed to be known, the update equation \eqref{eq: FP-update} is tractable since we omit the normalizing constant. However, it is sometimes beneficial to consider the normalized density as we discuss further in Section~\ref{section: normalizing constant}. It remains to solve the prediction equation \eqref{eq: FP-prediction}, also known as the Fokker--Planck equation or the Kolmogorov forward equation. In Section~\ref{Section: Methods} we recall two optimization-based methods that approximately solve the filtering equations, the deep splitting filter \cite{baagmark2024convergent} and the deep BSDE filter \cite{baagmark2025nonlinear}. The main objective of this work is comparing these methods that are constructed to mitigate the curse of dimensionality through Monte Carlo approximations and neural networks.

However, an intrinsic issue of solving for $p$ directly in \eqref{eq: FP-prediction}--\eqref{eq: FP-update}, is to handle the fact that probability density values diminish in increasing dimension. More precisely, with increasing state dimension $d$, the numerical values of $p$ decrease exponentially. This can be seen by considering an initial Gaussian distribution $\pi_0 = \N(0,I)$ in a $100$-dimensional space where the value of $\pi_0$ in the origin is approximately $10^{-40}$. Multiplying this with a likelihood of a similar structure quickly becomes unstable if not handled correctly. Moreover, since the methods that we consider rely on iterative gradient descent to optimize a set of parameters, they are further affected by the vanishing gradient phenomenon \cite{schoenholz2017deep}. Both of these issues can be mitigated by directly modeling the log-density of the filtering equations, which leads to positivity-preserving schemes for the density approximation. We state these equations next and prove them in Appendix~\ref{appendix: Proofs}.
\begin{theorem} \label{theorem: log equation}
    Let $p$ be the solution to \eqref{eq: FP-prediction}--\eqref{eq: FP-update}. Define the log-transformed solution by $v=(v_k)_{k=0}^K$, $v_k\colon [t_k,t_{k+1}]\times \R^d\times \mathbb{O}_k\to \R$, such that $v_k = -\log p_k$, for  $k=0,\dots,K$. The log-density $v$, satisfying $v_{0}(0)=-\log(\pi_0)$, is, for $k=0,\dots,K-1$, $x\in\R^d$, $o_{1:k}\in \mathbb{O}_k$, and $t\in[t_k,t_{k+1}]$, the solution to 
    \begin{align}
        \begin{split}
        \label{eq: log FP-prediction}
        &v_{k}(t,x,o_{1:k}) 
        = 
        v_{k}(t_{k},x,o_{1:k})
        \\
        &\hspace{5em}
        +
        \int_{t_{k}}^{t} 
        \Big(
        A v_{k}(s,x,o_{1:k}) 
        +
        f_{\mathrm{log}}
        (
        x,
        v_{k}(s,x,o_{1:k}),
        \nabla 
        v_{k}(s,x,o_{1:k})
        )
        \Big)
        \dd s,
        \end{split}
        \\
        \label{eq: log FP-update}
        &v_{k+1}
        (t_{k+1},x,o_{1:k+1})
        =
        v_{k}
        (t_{k+1},x,o_{1:k})
        -
        \log
        (
        L(o_{k+1},x)
        )
        ,
    \end{align}
    where $f_{\mathrm{log}}$ is defined, for $x\in\R^d$, $u\in\R$, and $w\in\R^d$, by
    \begin{align*}
        f_{\mathrm{log}}
        (
        x,
        u,
        w
        )
        &=
        -
        \frac{1}{2}
        \|
        \sigma(x)^\top
        w
        \|^2
        -
        f
        (
        x,
        1,
        -w
        ).
    \end{align*}
\end{theorem}
We note two aspects about \eqref{eq: log FP-prediction}--\eqref{eq: log FP-update}. The first and most important one is that the log-equation is also written in terms of the generator $A$, which means it is a type of nonlinear Kolmogorov backward equation. This makes both the deep splitting and the deep BSDE methods applicable. Secondly, the first order term $f_{\text{log}}$ is now independent of its second argument and nonlinear in its third argument. 

This formulation also has an important practical implication: approximating the solution in terms of the log-density naturally leads to positivity-preserving density approximations when the filter is reconstructed by exponentiation. In the next section, we introduce numerical approximation schemes for both the density and log-density formulations.

\section{Deep density methods}
\label{Section: Methods}
In this section we outline the methods under consideration, namely the deep splitting filter and the deep BSDE filter. In addition to these previously studied approaches, we also propose adaptations based on the logarithmic formulation \eqref{eq: log FP-prediction}--\eqref{eq: log FP-update}, which yield positivity-preserving schemes for the density approximation. Despite their different appearances, the derivations of these methods follow a very similar structure, and it is convenient to introduce a common probabilistic framework from which they can be developed. To this end, we define an auxiliary process $X \colon [0,T]\times \Omega \to \mathbb{R}^d$, that for all $t\in[0,T]$, $\mathbb{P}$-a.s., satisfies
\begin{align*} 
%\label{eq: auxiliary state}
    X_t 
    &=
    X_0 
    + 
    \int_0^t 
    \mu(X_s)
    \dd s 
    + 
    \int_0^t 
    \sigma(X_s)
    \dd W_s,
    \quad 
    X_0 \sim q_0.
\end{align*}
The particular choice of initial distribution \(q_0\) is not essential for the derivation, but a natural choice is to set \(q_0=\pi_0\), which ensures that $X$ equals $S$ in distribution. In some examples, however, we will find it advantageous to select \(q_0\) differently, for instance to increase the variance in the trajectories of $X$. Both methods are obtained by applying Itô’s formula to $p_k(t_k-t,X_t)$, which yields two closely related but distinct Feynman--Kac representations (corresponding to the two approximation strategies).

Since the methods are based on simulation, and closed-form solutions are only available in some specific cases, we approximate the process $X$. We fix a uniform time step $\tau = \frac{T}{KN} = t_{k,n+1}-t_{k,n}$, $k=0,\dots,K-1$, $n=0,\dots,N-1$, and approximate $X$, by $\mathcal{X}^k$, on each interval \([t_k,t_{k+1}]\) with the Euler--Maruyama scheme, for $n=0,\dots,N-1$, given by
\begin{align}
\label{eq: Euler--Maruyama SDE}
    \mathcal{X}_{n+1}^k
    &=
    \mathcal{X}_{n}^k
    +
    \mu(\mathcal{X}_{n}^k)\,(t_{k,n+1}-t_{k,n})
    +
    \sigma(\mathcal{X}_{n}^k)\,(W_{t_{k,n+1}}-W_{t_{k,n}}).
\end{align}
Successive intervals are linked by setting $\mathcal{X}_N^{k-1} = \mathcal{X}_0^k$, so that there is a one-to-one correspondence between the discretization times \(t_{k,n}\) and the indices \(\mathcal{X}_n^k\). This construction provides the numerical backbone on which both the deep splitting and deep BSDE filters are formulated. 

\subsection{Deep splitting filter}
\label{section: deep splitting}
The original deep splitting methodology was developed for parabolic PDEs in \cite{Arnulf_PDE}, analyzed in \cite{Germain}, later extended to a class of SPDEs in \cite{Arnulf} and to partial integro-diﬀerential equations in \cite{frey2022convergence}. In \cite{baagmark2024convergent} the method was applied to the filtering equations with suitable adaptations, and a strong convergence order of one was proven under a parabolic Hörmander condition with smooth and bounded coefficients. To remain coherent with the deep BSDE framework, we adopt a slightly modified, but equivalent, notation. The idea is to approximate the conditional predictive density through a recursive optimization problem. We define, for a fixed discretization time step $\tau$, the operator $G^\tau$ acting on a test function $\varphi \in C^1(\R^d;\R)$ by
\begin{align*}
    (G^\tau\varphi)(x)
    &=
    \varphi(x)
    +
    \tau
    f(x,
    \varphi(x),
    \nabla
    \varphi(x)
    ),
    \quad
    x\in\R^d
    .
\end{align*}
For every $k=0,\dots,K-1$ and $n=0,\dots,N-1$, the approximation $\phi_{k,n+1}$ is then defined, for $k=0,\dots,K-1$, $n=0,\dots,N-1$, as the solution of the recursive optimization problem 
\begin{align} 
\label{eq: deep splitting minimization}
\begin{split}
    &\mathop{\mathrm{min}}_{\phi \in 
    L^\infty(\mathbb{O}_k;
    C(\mathbb{R}^d;\mathbb{R}))}
    \mathbb{E}
    \bigg[ \Big| 
    \phi(\mathcal{X}_{n}^k,O_{1:{k}})
    -
    G^\tau
    \widetilde{g}_{k,n}
    (\mathcal{X}_{n+1}^k,O_{1:k})
    \Big|^2 \bigg],
\end{split}
\end{align}
where the function $\widetilde{g}_{k,n}$ encodes the recursion and is defined by
\begin{align*} 
    \widetilde{g}_{k,n}(x,o_{1:k})
    &=
    \begin{cases}
        \phi_{k,n}
        (x,o_{1:k})
        ,
        & 
        n\geq 1,
        \\
        \phi_{k-1,N}
        (x,o_{1:k-1})
        L(o_k,x)
        ,
        & 
        n=0,\ 
        k\geq 1,
        \\
        \pi_0(x),
        & 
        n=0,\ 
        k=0.
    \end{cases} 
\end{align*}
That is, for prediction steps we use the previous approximation $\phi_{k,n}$, while at observation times ($n=0$) we incorporate the likelihood $L$ through Bayes' formula as in \eqref{eq: FP-update}. This recursion yields an approximation of the conditional predictive density $p$ at the discretization times $t_{k,n}$, in the sense that $\phi_{k,n}(x,o_{1:k}) \approx p_k(t_{k,n},x,o_{1:k})$ for all $k=0,\dots,K-1$, $n=1,\dots,N$, $o_{1:k}\in \mathbb{O}_k$ and $x\in\R^d$. To this end, we define the approximative filter $\widetilde{p}= (\widetilde{p}_k)_{k=1}^K$, of the exact filter $(p_k(t_k))_{k=1}^K$, by
\begin{align*}
%\label{eq: deep splitting filter}
    \widetilde{p}_k
    (x,o_{1:k})
    &=
    \phi_{k-1,N}
    (x,o_{1:k-1})
    L(o_k,x)
    ,\quad
    o_{1:k}\in \mathbb{O}_k
    ,\
    x\in \R^d
    .
\end{align*}
In the logarithmic formulation of the method, the operator $G^\tau$ is modified to act on $\varphi \in C^1(\R^d;\R)$ by
\begin{align*}
    (G^\tau\varphi)(x)
    &=
    \varphi(x)
    +
    \tau
    f_{\mathrm{log}}
    (x,
    \varphi(x),
    \nabla
    \varphi(x)
    ),
    \quad
    x\in\R^d
    ,
\end{align*}
while the recursion $\widetilde{g}_{k,n}$, for $x\in\R^d$ and $o_{1:k}\in\mathbb{O}_k$, becomes
\begin{align*} 
    \widetilde{g}_{k,n}(x,o_{1:k})
    &=
    \begin{cases}
        \phi_{k,n}
        (x,o_{1:k})
        ,
        & 
        n\geq 1,
        \\
        \phi_{k-1,N}
        (x,o_{1:k-1})
        -
        \log
        (L(o_k,x))
        ,
        & 
        n=0,\ 
        k\geq 1,
        \\
        -
        \log
        (
        \pi_0(x)
        ),
        & 
        n=0,\ 
        k=0.
    \end{cases} 
\end{align*}
We represent $\phi$ in \eqref{eq: deep splitting minimization} by a neural network, approximate the expectation using Monte Carlo samples of $(\mathcal{X}, O)$, and employ stochastic gradient descent to locate the minimizer. Details on the network architectures and optimization settings are provided in Supplementary~Materials~\ref{SM: Training}.

In summary, the deep splitting method reformulates the filtering equations into a sequence of recursive optimization problems, where each step approximates the action of the operator $G^\tau$ on the predictive densities. The recursion alternates between propagation by the Euler--Maruyama scheme and correction through Bayes' formula at observation times. In this way, the approximations $(\phi_{k,n})_{k=0,n=1}^{K-1,N}$ provide a discretized representation of the conditional predictive densities, while the terminal quantities $(\widetilde{p}_k)_{k=1}^K$ yield an explicit approximation of the (unnormalized) filter itself.

\subsection{Deep BSDE filter}
\label{section: deep BSDE}
In contrast to the deep splitting scheme outlined above, the deep BSDE formulation leads to a single optimization problem per observation time. The method was introduced in \cite{E_2017} for high-dimensional semilinear PDEs, and has been extended in various directions \cite{andersson2026deep,beck2019machine, chan2019machine,han_convergence}; see the survey \cite{han2025brief} for an overview. In the context of nonlinear filtering, the approach was adapted in \cite{baagmark2025nonlinear} by coupling deep BSDE formulations for the prediction step with Bayes’ formula for the update, resulting in an iterative alternation between propagation and correction. The approximation error was bounded, under an Hörmander condition, by an estimate that combines an a priori $1/2$-order discretization error with an a posteriori term given by the attained objective value.

For each $k=0,\dots,K-1$, the optimization problem is designed such that its solution $\phi_k$ approximates the conditional predictive density $p_k(t_{k+1})$. We solve
\begin{align} \label{eq: Deep BSDE minimization}
\begin{split}
    &\hspace{-4em}
    \underset{
    \substack{
    \phi\in 
    L^\infty(\mathbb{O}_k;
    C(\R^d
    ;\R))
    \\
    (\overline{v}_n)_{n=0}^{N-1}\in 
    L^\infty(\mathbb{O}_k;
    C(\R^d
    ;\R^d))^N
    }}
    {\mathrm{min}}\  
    \E
    \Bigg[
    \Big|
    \mathcal{Y}_{k,N}^{O_{1:k}}
    -
    \overline{g}_k
    (\mathcal{X}_{N}^k,
    O_{1:k})
    \Big|^2
    \Bigg],
\end{split}
\end{align}
where the process $(\mathcal{Y}_{n}^{O_{1:k}})_{n=0}^N$ is defined recursively, with $n=0,\dots,N-1$, by an Euler--Maruyama discretization of the BSDE
\begin{align*}
    \mathcal{Y}_{0}^{O_{1:k}}
    &=
    \phi(\mathcal{X}_{0}^k, O_{1:k}),
    \\
    \begin{split}
        \mathcal{Y}_{n+1}^{O_{1:k}}
        &=
        \mathcal{Y}_{n}^{O_{1:k}}
        -
        f(
        \mathcal{X}_{n}^k,
        \mathcal{Y}_{n}^{O_{1:k}},
        \overline{v}_{n}
        (\mathcal{X}_{n}^k, O_{1:k})
        )
        (t_{k,n+1}-t_{k,n})
        \\
        &\hspace{6em}+
        \overline{v}_{n}
        (\mathcal{X}_{n}^k, O_{1:k})^\top
        \sigma
        (\mathcal{X}_{n}^k)
        (W_{t_{k,n+1}}-W_{t_{k,n}})
        .
    \end{split}
\end{align*}
Here the auxiliary functions $(\overline v_\ell)_{\ell=0}^{N-1}$, multiplied with $\sigma$, act as the approximations to the $Z$-components in the continuous-time BSDE formulation. The terminal condition $\overline{g}_k$ incorporates the update step via Bayes' formula, and is defined from the previous optimum $\phi_{k-1}$, for $x\in\R^d$ and $o_{1:k}\in\mathbb{O}_k$, by
\begin{align*} 
%\label{eq: overline-g_k}
    \overline{g}_k(x,o_{1:k})
    &=
    \begin{cases}
        \phi_{k-1}
        (x,o_{1:k-1})
        L(o_k,x),
        & 
        k\geq 1,
        \\
        \pi_0(x),
        & 
        k=0.
    \end{cases} 
\end{align*}
We define the deep BSDE filter approximation $\widehat{p} = (\widehat{p}_k)_{k=1}^K$, of the exact filter $(p_k(t_k))_{k=1}^K$, by
\begin{align*} 
%\label{eq: deep bsde filter}
    \widehat{p}_k
    (x,o_{1:k})
    &=
    \phi_{k-1}
    (x,o_{1:k-1})
    L(o_k,x)
    ,\quad
    o_{1:k}\in \mathbb{O}_k
    ,\
    x\in \R^d
    ,
\end{align*}
where $(\phi_k)_{k=0}^{K-1}$ are the obtained solutions of the optimization problems \eqref{eq: Deep BSDE minimization}. Similar to the deep splitting approach, we also outline the version approximating the log-density \eqref{eq: log FP-prediction}--\eqref{eq: log FP-update}. In the logarithmic version we switch $f$ for $f_{\mathrm{log}}$ and define $\overline{g}_k$, for $x\in\R^d$ and $o_{1:k}\in\mathbb{O}_k$, by
\begin{align*} 
    \overline{g}_k(x,o_{1:k})
    &=
    \begin{cases}
        \phi_{k-1}
        (x,o_{1:k-1})
        -
        \log
        (
        L(o_k,x)
        ),
        & 
        k\geq 1,
        \\
        -
        \log
        (
        \pi_0(x)
        ),
        & 
        k=0.
    \end{cases} 
\end{align*}
As in the deep splitting method, $\big(\phi,(\overline{v}_n)_{n=0}^{N-1}\big)$ in \eqref{eq: Deep BSDE minimization} are represented by deep neural networks, the expectation is estimated by Monte Carlo samples of $(\mathcal{X}, O)$, and the minimizer is found using stochastic gradient descent. Further details on architectures and optimization are given in Supplementary~Materials~\ref{SM: Training}.

In summary, the deep BSDE method approximates the prediction step by discretizing a backward stochastic differential equation associated with the nonlinear filtering PDE. The optimization simultaneously learns the predictive density $\phi_k$ and the auxiliary processes $(\overline v_n)$, while the update step is enforced through the terminal condition $\overline g_k$. Thus the method condenses the filtering problem into one global optimization per observation time, in contrast to the stepwise recursion of the deep splitting scheme.

\section{Experimental setup}
\label{section: experiment setup}
In this section we briefly outline and summarize details regarding implementation and evaluation of the experiments that we report in Section~\ref{Section: Experiments}.

\subsection{Normalizing constant}
\label{section: normalizing constant}
The filtering equations \eqref{eq: FP-prediction}--\eqref{eq: FP-update}, and likewise \eqref{eq: log FP-prediction}--\eqref{eq: log FP-update}, yield solutions in the form of unnormalized densities (or log versions thereof). Theoretically, and in the corresponding analysis of the methods \cite{baagmark2025nonlinear,baagmark2024convergent}, this unnormalized setting is sufficient.
However, for large numbers of observation times $K$, it is often beneficial to apply an approximate normalization to enhance numerical stability.
This occurs in the update step \eqref{eq: FP-update} when multiplying the predicted density with the likelihood, often with very small overlap of probability mass between the two, resulting in a much smaller total probability mass. Without normalization, the resulting density deteriorates numerically when the procedure is repeated many times.

In one dimension it is straightforward to perform a quadrature approximation of the normalization constant. However, as we are mostly interested in scaling to higher dimensions, we tackle the issue by importance sampling. A natural choice of importance distribution is a variance-scaled extended Kalman filter; in addition, we also consider a cheaper observation-independent wide-tailed Gaussian proposal. Even though in many examples we expect the approximate Kalman filter to perform poorly, it still gives an idea of the domain of interest, and increasing the variance yields effective samples in the domain of the true filter. See Supplementary~Materials~\ref{SM: Normalization} for details.

\subsection{Training}
During training, the aim is to approximately solve the optimization problems \eqref{eq: deep splitting minimization} and \eqref{eq: Deep BSDE minimization} for  DSF and BSDEF respectively. If nothing else is stated we have used Fully Connected neural Networks (FCN) to parameterize the networks in the optimization problems. We employ the deep approximative methods for both the original formulation and the log-formulation, yielding four methods of interest, denoted DSF, BSDEF, LogDSF, and LogBSDEF. The details about architectures and training can be found in Supplementary~Materials~\ref{SM: Training}, which includes hyperparameters such as the number of discretizations steps $N$ that is used.

\subsection{Classical methods}
In Section~\ref{Section: Experiments} we benchmark against three classical baselines: the extended Kalman filter, the ensemble Kalman filter, and the bootstrap particle filter. The EKF linearizes the nonlinear model using a first-order Taylor expansion and then applies the Kalman filter, yielding a computationally cheap baseline. The EnKF and PF propagate ensembles/particles forward; achieving a fixed accuracy typically requires an exponentially growing number of ensemble members/particles as the state dimension increases, making these methods computationally expensive. See Supplementary~Materials~\ref{SM: Classical methods} for details on the three methods.

\subsection{Error metrics and evaluation}
To capture the differences in the performance of the methods, we investigate moment-based errors and uncertainty-based errors. The first type describes how well the methods manage to capture the mean or how well the computed mean estimates the true trajectory. The second type describes how well we capture the uncertainty contained in the overall density, including the tails of the distribution. 

The first moment-based metric is a First Moment Error (FME), which measures the error between the true conditional mean, associated with the solution to \eqref{eq: FP-prediction}--\eqref{eq: FP-update}, and the approximate mean. The second moment-based metric is a Mean Absolute Error (MAE), which measures the difference between the true trajectory $S$ and the approximative mean. More precisely, we simulate $M$ independent trajectories $(S_{t_k}^{(m)},O_k^{(m)})_{k=1}^K$, $m=1,\dots,M$, from \eqref{eq: state}--\eqref{eq: obs}. For each $O_k^{(m)}$ we evaluate the true conditional mean $\mu_k^{(m)}$, as well as an approximate mean $\widehat{\mu}_k^{(m)}$ from each considered approximate filter. We define the metrics, for each $k=1,\dots,K$, and $\widehat{\mu}_k$, by
\begin{align*}
    \text{FME}(\widehat{\mu}_k)
    &=
    \mathbb{E}_{O_k}
    \Big[
    \big\|
    \mu_k
    -
    \widehat{\mu}_k
    \big\|
    \Big]
    \approx
    \frac{1}{M}
    \sum_{m=1}^{M}
    \big\|
    \mu_k^{(m)}
    -
    \widehat{\mu}_k^{(m)}
    \big\|
    ,
    \\
    \text{MAE}(\widehat{\mu}_k)
    &=
    \mathbb{E}_{S_{t_k},O_k}
    \Big[
    \big\|
    S_{t_k}
    -
    \widehat{\mu}_k
    \big\|
    \Big]
    \approx
    \frac{1}{M}
    \sum_{m=1}^{M}
    \big\|
    S_{t_k}^{(m)}
    -
    \widehat{\mu}_k^{(m)}
    \big\|
    .
\end{align*}
While $\text{FME}(\mu_k)=0$ by definition, $\text{MAE}(\mu_k) \neq 0$ remains useful as a scale: the mean of the true filter, $\mu_k$, is the $L^2(\Omega,\mathfrak{S}(O_{1:k});\mathbb{R}^d)$-optimal estimator (the conditional expectation), so its MAE provides the natural reference level and lower bound. We therefore define the relative MAE (rMAE),
\begin{align*}
    \text{rMAE}(\widehat{\mu}_k)
    &=
    \frac{\text{MAE}(\widehat{\mu}_k) - \text{MAE}(\mu_k)}
    {\text{MAE}(\mu_k)},
\end{align*}
which equals $0$ for the optimal estimator and is smaller-is-better in general. We report this metric in percentage to highlight how much larger the error of the approximation $\widehat{\mu}_k$ is compared to the optimal estimator $\mu_k$. 

For the uncertainty-based metrics, which quantify distributional discrepancy, we use the forward Kullback--Leibler divergence when a reliable reference is available. Recall that $p_k$ solves \eqref{eq: FP-prediction}--\eqref{eq: FP-update}, and let $\widehat{p}_k$ denote the density of an approximation. For each $k=1,\dots,K$, and $\widehat{p}_k$, we define
\begin{align*}
\begin{split}
    \text{KLD}(\widehat{p}_k)
    &=
    \mathbb{E}_{O_{1:k}}
    \Big[
    D_{\text{KL}}
    \big(
    p_k
    (t_k,\cdot,O_{1:k})
    \,
    \big\|
    \,
    \widehat{p}_k
    (\cdot,O_{1:k})
    \big)
    \Big]
    \\
    &=
    \mathbb{E}_{O_{1:k}}
    \Big[
    \int_{\R^d}
    \log
    \bigg(
    \frac{p_k(t_k,x,O_{1:k})}
    {\widehat{p}_k(x,O_{1:k})}
    \bigg)
    p_k(t_k,x,O_{1:k})
    \dd x
    \Big].
\end{split}
\end{align*}
In practice, the approximation of the Kullback--Leibler divergence relies on approximating the integral over $\R^d$, with $J$ Monte Carlo samples $(Z^{(j,m)})_{j=1}^J$ from the reference distribution $p_k(t_k,x,O_{1:k}^{(m)})$ for every $m=1,\dots,M$, according to
\begin{align*}
\begin{split}
    \text{KLD}(\widehat{p}_k) 
    \approx
    \frac{1}{M}
    \sum_{m=1}^{M}
    \frac{1}{J}
    \sum_{j=1}^J
    \log
    \Bigg(
    \frac{
    p_k(t_k,Z^{(j,m)},O_{1:k}^{(m)})
    }
    {
    \widehat{p}_k(Z^{(j,m)},O_{1:k}^{(m)})
    }
    \Bigg)
    .    
\end{split}
\end{align*}
In the linear-Gaussian setting, where the Kalman filter provides an analytical solution, this is straightforward, since we can sample $Z$ directly from the Gaussian posterior $p_k$. In the low-dimensional nonlinear setting, where we use a bootstrap particle filter as the reference, we approximate sampling from $p_k$ by drawing $Z$ from the empirical particle distribution at time $t_k$, using the normalized particle weights as sampling probabilities. This procedure ensures that the Monte Carlo estimate of the Kullback--Leibler divergence remains consistent with the reference solution in both the linear and low-dimensional nonlinear cases.

However, obtaining a reliable reference solution for the full density is generally infeasible in the high-dimensional nonlinear setting. In such cases, instead we benchmark the uncertainty-based performance of an approximation through the Negative Log-Likelihood (NLL) of the ground truth state. The NLL measures how much probability mass the filter assigns to the realized state, and its expectation corresponds to the cross-entropy between the true and approximate posteriors, differing from the KL divergence only by an additive constant \cite{CoverThomas2006}. Hence, smaller NLL values indicate better alignment with the true filtering distribution. This metric is not zero for the optimal filter, but we rather seek as small values as possible. The metric is defined, for each $k=1,\dots,K$, and $\widehat{p}_k$, by
\begin{align*}
\begin{split}
    \text{NLL}
    (\widehat{p}_k)
    &=
    \mathbb{E}_{S_{t_k},O_{1:k}}
    \Big[
        -
        \log\big(
            \widehat{p}_k
            (S_{t_k}
            , 
            O_{1:k})
        \big)
    \Big]
    \approx
    -
    \frac{1}{M}
    \sum_{m=1}^{M}
    \log\big(
        \widehat{p}_k
        (S_{t_k}^{(m)}
        ,
        O_{1:k}^{(m)})
    \big).
\end{split}
\end{align*}
In Supplementary~Materials~\ref{SM: Evaluation} we report all parameters, normalization methods and references used for each example from Section~\ref{Section: Experiments}.

\section{Experiments}
\label{Section: Experiments}
In this section we benchmark the proposed deep filters, DSF and BSDEF, and their logarithmic variants across problems of increasing difficulty, using the methodology of Section~\ref{Section: Methods} and comparing to PF, EKF, and EnKF. 
While convergence results, as $N\to\infty$, are available for the deep density methods under strong regularity conditions, these conditions are often not met in the examples considered here. Nonetheless, the methodology of Section~\ref{Section: Methods} yields well-defined algorithms, and we use the experiments to quantify their empirical performance and robustness compared to the classical approaches.
We begin, in Section~\ref{section: toy}, with two one-dimensional toy problems with reliable references: a linear Ornstein--Uhlenbeck process and a nonlinear bistable system. In Section~\ref{section: longrange highdim OU}, robustness and applicability in linear high-dimensional settings are assessed through a long-horizon Ornstein--Uhlenbeck process and a high-dimensional counterpart, both serving as simple canonical tests. 
Section~\ref{section: schlogl} continues with the Schl\"ogl model providing a first step toward physically grounded dynamics. In Section~\ref{section: L96} we test the models on a stochastic Lorenz-96 system, using high-accuracy PF surrogates when feasible. Performance is reported over time via $\mathrm{MAE}$, $\mathrm{FME}$, and $\mathrm{KLD}$, with particular attention to the stability and accuracy gains of the logarithmic formulations. The final Section~\ref{section: computional efficiency} analyzes computational scaling and inference-time costs. Normalization, architectural, training, and evaluation details appear in Supplementary~Materials~\ref{SM: Normalization}--\ref{SM: Evaluation}.

We systematically compare deep density based methods with classical approaches by varying the particle count in PFs and the ensemble size in EnKFs. To keep the presentation focused, we report configurations that highlight performance differences and omit less informative variants.

\subsection{Toy examples}
\label{section: toy}
We consider two simple one-dimensional equations, one linear and one nonlinear. In both of the examples we set $T=1$, $K=10$, $h(x) = x$, $R=1$, $\sigma(x) = 1$ and $\pi_0 = \N(0,1)$. The linear equation has a mean reverting drift $\mu(x) = -x$ and the corresponding SDE is solved by the Ornstein--Uhlenbeck process. In this case, the reference solution to \eqref{eq: FP-prediction}--\eqref{eq: FP-update} is provided by the Kalman filter. The drift in the nonlinear equation, $\mu(x) = -\frac{2}{5}(5x-x^3)$, is the negative gradient of a double-well potential. Hence the unconditional distribution of the solutions becomes bimodal and we refer to this process as a bistable process. This example, but with discrete dynamics, was previously studied in \cite{bao2024score}. In this case we use a particle filter with $10^6$ particles and $128$ intermediate time steps between observations to sufficiently well approximate the true filter. 

\begin{figure}[h]
    \centering
    \begin{minipage}[t]{0.9\linewidth}  
    \raggedright
    
    \input{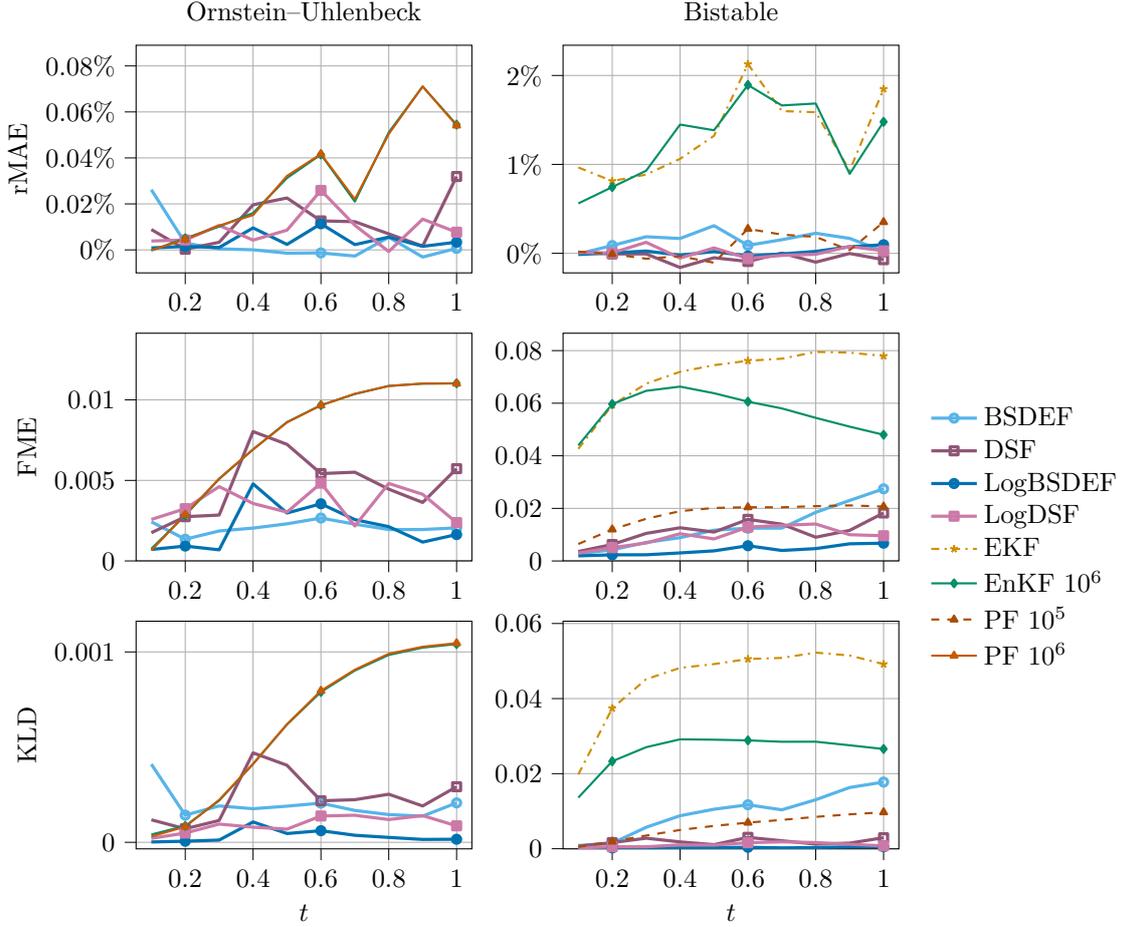}
    \end{minipage}%
    \begin{minipage}[t]{0.1\linewidth} 
    \vspace{0pt}
    \raggedright
    \begin{tikzpicture}[overlay]
    \node at (0,0) {$ $};
    \node at (-5.0,-0.5) {\ref*{legend: OU MAE}};
    \end{tikzpicture}
    \end{minipage}
    \vspace{2em}
    \caption{On the left and right panels the results for the Ornstein--Uhlenbeck process and the bistable process are depicted, respectively. From top to bottom the rMAE, FME, and KLD metrics are illustrated.}
    %\vspace{-2em}
    \label{fig: toy}
\end{figure}

In Figure~\ref{fig: toy} we see the MAE, FME and KLD evaluated on the two examples. We compare the performance to particle filters with a range between $10^4$ and $10^6$ particles with one intermediate time step, and in the nonlinear case also to the EKF. Both the DSF and the BSDEF perform well but we note that the logarithmic versions, particularly the LogBSDEF, achieve a lower error in almost all metrics. It should be stressed here that the computational time for estimation and density evaluation are about two and five orders of magnitudes faster, respectively, see Section~\ref{section: computional efficiency}.

\subsection{High-dimensional Ornstein--Uhlenbeck}
\label{section: longrange highdim OU}
To bridge the toy problems and large-scale tests, we study a controlled linear benchmark that isolates temporal and dimensional effects. The Ornstein--Uhlenbeck process serves this purpose, since it admits an analytic Kalman reference and its unconditional density is close to stationarity at $t=0$ for $\pi_0 = \N(0,I)$. We consider a long-horizon example with $T=10$, $K=100$, and $d=10$, and a high-dimensional example with $T=1$, $K=10$, and $d=100$. In both examples we have $d'=d$, $h(x)=x$, $R=I$, $\mu(x)=-x$, $\sigma(x)=I$, and $\pi_0=\N(0,I)$. Furthermore, in the long-horizon example we compare two different architectures, namely in addition to the standard FCN, we also employ a Long Short-Term Memory (LSTM)-adaptation, for the LogBSDEF, suitable for longe-horizon problems. See Supplementary~Materials~\ref{SM: Architectures} for details on the architectures. 

\begin{figure}
    \centering
    \begin{minipage}[t]{0.9\linewidth}  
    \raggedright
    
    \input{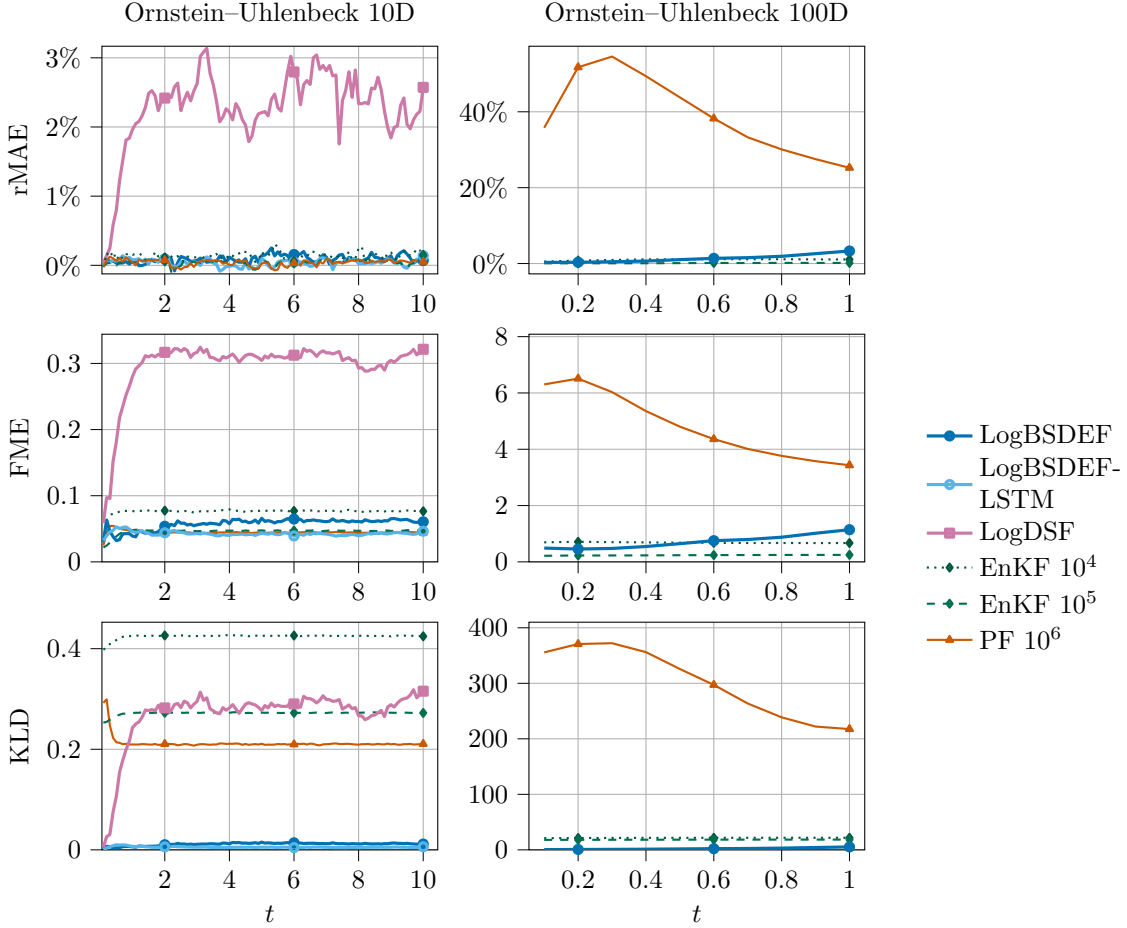}
    \end{minipage}%
    \begin{minipage}[t]{0.1\linewidth}  
    \vspace{0pt}
    \raggedright
    \begin{tikzpicture}[overlay]
    \node at (0,0) {$ $};
    \node at (-5.0,-0.6) {\ref*{legend: OU LR MAE group}};
    \end{tikzpicture}
    \end{minipage}
    \vspace{2em}
    \caption{On the left and right panels the results for the long-horizon $10$-dimensional Ornstein--Uhlenbeck process and the short-horizon $100$-dimensional Ornstein--Uhlenbeck process are depicted, respectively. From top to bottom the rMAE, FME, and KLD metrics are illustrated.}
    \vspace{-1.0em}
    \label{fig: OU long high}
\end{figure}

In Figure \ref{fig: OU long high} we evaluate the deep learning based methods as well as PFs and EnKFs with different number of particles. Looking at the long-horizon process we see, in particular, that the two LogBSDEF models perform at least as good as the particle filter for the rMAE and FME, and for the KLD it outperforms the classical approximations by a large margin. Moreover, the approximation quality remains stable over time and by introducing the LSTM architecture we make LogBSDEF perform best out of all benchmarked methods. In Supplementary~Materials~\ref{SM: densities} we demonstrate the methods on a single trajectory by plotting trajectories and marginal densities, where we also see a satisfactory performance of the LogDSF. 

For the 100-dimensional Ornstein--Uhlenbeck example we clearly see that the particles filters give quite poor approximations even when using $10^6$ particles. The LogBSDEF gives satisfactory results, where the averaged performance is on par with the ensemble Kalman filters, which are expected to yield excellent performance in this linear setting. Remark that we only include the results for LogDSF and LogBSDEF, since the non-logarithmic variants are unstable during training due to the previously mentioned numerical problems with high-dimensional probability densities. Similarly, we found that the LogDSF is unstable in the $100$-dimensional example and thus omit it there. In the figures below, whenever a method is excluded, it means that it failed to train properly. 

The computational gain from using the LogBSDEF compared to the underperforming PF with $10^6$ particles, is even higher than for the one-dimensional example, namely three orders of magnitude faster, see Section~\ref{section: computional efficiency}. From this point on we refrain from commenting on more results on computational efficiency. It depends on dimension and not the problem at hand. 

\subsection{Schl\"ogl model}
\label{section: schlogl}
The nonlinear SDE considered in this example originates from the classical Schlögl model \cite{schlogl1972chemical, vellela2009stochastic, Vlysidis2018}, a canonical benchmark for bistable reaction systems. The model describes the concentration process $S$ of a single reacting species subject to autocatalytic reactions and exchange with two external reservoirs. Using the chemical Langevin approximation \cite{gillespie2000chemical}, one can derive the dynamics satisfying \eqref{eq: state} with $\sigma\colon \R \to \R^{1\times 4}$ and a four-dimensional Brownian motion $B = (B^{(i)})_{i=1}^4$ corresponding to the four reaction channels. The dynamics are governed by the drift and diffusion coefficients
\begin{align*}
    \mu(x) 
    &= 
    \big(
    b_1(x)
    - b_2(x) 
    + b_3(x) 
    + b_4(x)
    \big)
    ,
    \\
    \sigma(x) 
    &= 
    \Big(
    \sqrt{b_1(x)},
    -\sqrt{b_2(x)},
    \sqrt{b_3(x)},
    -\sqrt{b_4(x)}
    \Big),
\end{align*}
where 
\begin{align*}
    b_1(x) 
    =
    \tfrac{\theta_1 A}{2} 
    x (x - 1)
    ,
    \quad
    b_2(x) 
    =
    \tfrac{\theta_2}{6}
    x (x - 1)(x - 2)
    ,
    \quad
    b_3(x) 
    =
    \theta_3 B
    ,
    \quad
    b_4(x) 
    =
    \theta_4 x
    .
\end{align*}
The reaction rate parameters are fixed to $\theta =(3\times 10^{-7},10^{-4},10^{-3},3.5)$, $A = 10^{5}$, and $B = 2\times 10^{5}$.
This setup produces a bistable system with two distinct equilibrium states separated by a potential barrier. In the Euler--Maruyama discretization \eqref{eq: Euler--Maruyama SDE}, the four Brownian channels can be merged into a single effective noise term, producing a scheme, identical in distribution, that reduces the computational cost. A more detailed presentation of the model is given in Supplementary~Materials~\ref{SM: schlogl chemistry}.

We fix $d' = 1$, $T=1$, $K=10$, $\pi_0 = \frac{\N(150,10^2) + \N(350,10^2)}{2}$, $h(x) = \log(1+x)$, and $R=0.5$ for the experiment. In the derivation of the deep filtering methods of Section~\ref{Section: Methods} we introduced $X$ to be initialized from a density $q_0$. In this example we set $q_0$ different from $\pi_0$ to account for the strong drift of the process. More precisely, we set $q_0 = \frac{\N(150,25^2) + \N(375,60^2)}{2}$ to better match the generated samples of $X_0$ with the distribution of $S_{t_1}$. See \cite[Section~3.4]{bagmark_1} for an extensive explanation of a very similar problem of matching densities. We use a particle filter with $10^6$ particles and $128$ intermediate time steps as a reference solution. 

In the experiment we employ both the LogBSDEF and the LogDSF, but the latter does not converge for the architectures and discretizations that we used, clearly indicating the difficulty of this problem. In the left column of Figure~\ref{fig: nonlinear metrics} we show the metrics over time, where we compare the LogBSDEF method to particle filters, ensemble Kalman filters, and the extended Kalman filter. Clearly the nonlinearity is too strong for the EKF and EnKFs to yield sufficiently good results. Looking at the FME, our proposed LogBSDEF method performs about as well as the particle filter with $10^5$ particles initially but the error accumulates over time. Similar to the other examples our method excels with respect to the KLD, where performance is mostly on par with the particle filter with $10^5$ particles.

% \begin{figure}[h]
% \captionsetup{justification=raggedleft}
% \input{include/figures/schlogl/groupplot}
% \begin{tikzpicture}[overlay]
% \node at (0,0) {$ $};
% \node at (-2.5,3.0) {\ref*{legend: schlogl MAE group}};
% \end{tikzpicture}
% \captionsetup{justification=justified}
% \vspace{10pt}
% \caption{Metrics for the Schl\"ogl model, shown left to right: rMAE, FME, and KLD.}
% \label{fig: schlogl metrics}
% \end{figure}

\subsection{Lorenz-96}
\label{section: L96}
In our final example, we tackle a strongly nonlinear, high-dimensional system, precisely the regime where classical methods succumb to the curse of dimensionality. The Lorenz-96 model is a high-dimensional chaotic dynamical system originally introduced in \cite{lorenz1996predictability} as a testbed for numerical weather prediction. It captures essential features of spatio-temporal chaos, while remaining simple enough for controlled experiments. In its deterministic form, the $d$-dimensional Lorenz-96 system satisfies \eqref{eq: state} with diffusion coefficient $\sigma(x) = 0$ and drift coefficient
\begin{align*}
    \mu(x)_i 
    &=
    (x_{i+1}-x_{i-2})
    x_{i-1}
    -
    x_i
    +
    F
    ,\quad
    i=1,\dots,d,
\end{align*}
where the indices are taken modulo $d$, that is $x_i = x_{i+d}$. The parameter $F>0$ acts as a constant external forcing, typically chosen in the range $F \in [4,16]$ to produce chaotic dynamics. To account for unresolved processes or model error, we consider a stochastic extension of the Lorenz-96 model by adding additive noise in the form of a Brownian motion with $\sigma(x) = \sigma I$, with $\sigma\in\R$.

The Lorenz-96 SDE serves as a standard testbed for assessing data assimilation techniques, due to its balance of chaotic behavior, high dimensionality, and computational tractability \cite{wilks2005effects}. It is extensively studied in the context of spatio-temporal chaos (e.g., fractal dimension scaling) \cite{karimi2010extensive}, bifurcation structure and wave propagation \cite{van2019symmetries}, as well as various data assimilation and neural network hybrid methods \cite{brajard2020combining}. 

We fix $d = 4$, $d' = 4$, $h(x)=x$, $T=1$, $K=10$, $\pi_0 = \N\big((F,F,F,F), I\big)$, $\sigma = I$, $R=2I$ and $F=8$. In this strongly nonlinear setting we do not have a good reference solution for increasing dimensions, but in the case of $d=4$, we found the particle filter with $10^6$ particles enhanced with $128$ intermediate time steps to be sufficiently good. We tried the LogBSDEF and the LogDSF, and compared their performance to the PF, EKF, and EnKF. In the experiments we tried a multitude of settings to make the LogDSF stable, but it consistently failed to converge. In Supplementary~Materials~\ref{SM: Training} we detail the settings we used for the LogBSDEF. 

In the right column of Figure~\ref{fig: nonlinear metrics} we demonstrate the results of the rMAE, FME, and KLD. The standard approximative particle filter diverges after some timesteps, where especially the KLD explodes. In the evaluation of the KLD, we clipped the values at $|\log(10^{-200})| \approx 460$, for numerical reasons, which is seen by the flattened curve of the PF. We see that our method performs superior to the EnKF with $10^6$ ensemble members in all metrics, but we can also note that none of the methods reach close to the MAE-minimum that the reference solutions obtain, indicating the severe difficulty of the problem. In Supplementary~Materials~\ref{SM: densities} we illustrate for a single trajectory, the state $S$ and corresponding filter means, projected on the $x_1x_2$- and $x_3x_4$-planes, respectively, to give an idea about the chaoticicity of the system.

\begin{figure}[h]
\captionsetup{justification=raggedleft}
\input{include/figures/nonlinear/groupplot}
\begin{tikzpicture}[overlay]
\node at (0,0) {$ $};
\node at (-6,-0.5) {\ref*{legend: schlogl MAE group}};
\end{tikzpicture}
\captionsetup{justification=justified}
\vspace{1.5em}
\caption{
On the left and right panels the results for the Schl\"ogl model and the $4$-dimensional Lorenz-96 model are depicted, respectively. From top to bottom the rMAE, FME, and KLD metrics are illustrated.}
\vspace{-2em}
\label{fig: nonlinear metrics}
\end{figure}

We continue by increasing the state dimension to $d=[4,10,20,40,100]$. In addition, we only have partial observations with $d'=[4,5,5,10,25]$, where we observe every, every second or every fourth position, respectively. More precisely, the measurement function is defined by 
\begin{align*}
    \big(h(x)\big)_i
    &=
    \begin{cases}
        &x_{i},\ \quad (d,d')=(4,4),\\
        &x_{2i},\quad (d,d')=(10,5),\\
        &x_{4i},\quad (d,d')=(20,5),\ (40,10),\ (100,25),
    \end{cases}
\end{align*}
for $i=1,\dots,d'$. The remaining parameters are the same as in the $d=4$ case. Here it is challenging to obtain reference solutions and instead we measure how the error increases by dimension in both the MAE and NLL metrics. Optimally, the performance is stable with increasing state dimension. We accumulate the scores for all the time steps for each dimension in both metrics, $\text{MAE} = \frac{1}{K}\sum_{k=1}^K \text{MAE}(\widehat{p}_k)$, $\text{NLL} = \frac{1}{K}\sum_{k=1}^K \text{NLL}(\widehat{p}_k)$ and report the metrics in Figure~\ref{fig: Lorenz highdim} on log-log axes. Since the NLL grows approximately linearly with the state dimension $d$ and the MAE with $d^{\frac{1}{2}}$, we plot two references slopes for each metric. 
In this example we ideally wanted to use more particles for the PFs and ensembles for the EnKFs, but we were limited by memory consumption. The classical methods scale poorly in this regard, especially when recovering the density function, see more on the computational efficiency in Section~\ref{section: computional efficiency}.
We can see that the PF underperforms severely in almost all dimensions. Furthermore, the performance of the LogBSDEF and the EnKF is mostly on par for up to about $d=20$, and then the EnKF starts diverging as well. The robust performance of the LogBSDEF demonstrates the potential of deep density methods to extend nonlinear filtering to high-dimensional settings.

\begin{figure}[h]
    \centering
    % This file was created with tikzplotlib v0.10.1.
\begin{tikzpicture}

\definecolor{crimson2143940}{RGB}{214,39,40}
\definecolor{darkgray176}{RGB}{176,176,176}
\definecolor{darkorange25512714}{RGB}{255,127,14}
\definecolor{darkturquoise23190207}{RGB}{23,190,207}
\definecolor{forestgreen4416044}{RGB}{44,160,44}
\definecolor{goldenrod18818934}{RGB}{188,189,34}
\definecolor{gray127}{RGB}{127,127,127}
\definecolor{lightgray204}{RGB}{204,204,204}
\definecolor{mediumpurple148103189}{RGB}{148,103,189}
\definecolor{orchid227119194}{RGB}{227,119,194}
\definecolor{sienna1408675}{RGB}{140,86,75}
\definecolor{steelblue31119180}{RGB}{31,119,180}

\begin{groupplot}[group style={group size=2 by 1,horizontal sep=1.5cm,
    vertical sep=1.5cm}, width=7.0cm,
  height=5.0cm]
\nextgroupplot[
  legend cell align={left},
  legend style={
    fill opacity=0.8,
    draw opacity=1,
    text opacity=1,
    at={(0.03,0.97)},
    anchor=north west,
    draw=none,
    nodes={scale=1.0, transform shape},
    legend columns=-1
  },
  tick align=outside,
  tick pos=left,
  title={MAE},
  x grid style={darkgray176},
  xlabel={$d$},
  xmajorgrids,
  xmin=3, xmax=110,
  xmode=log,
  log basis x=2,
  log ticks with fixed point,
  xtick={4,10,20,40,100},
  xtick style={color=black},
  y grid style={darkgray176},
  ymajorgrids,
  ymode=log,
  log basis y=10,
  log ticks with fixed point,
  ymin=1e0, ymax=2e3,
  ytick style={color=black},
  legend to name={legend: Lorenz highdim MAE group},
]

% --- data (unchanged) ---
% \addplot [semithick, steelblue31119180, mark=*, mark size=1, mark options={solid}]
\addplot[method LogBSDEF, with marks every]
table {%
4   3.3649740222569315
10  15.003372873273573
20  22.169335354248496
40  31.701279888115796
100 50.57628660082919
};
\addlegendentry{LogBSDEF}

% \addplot [semithick, crimson2143940, mark=*, mark size=1, mark options={solid}]
\addplot[method EnKF 1e6, with marks every]
table {%
4   4.367744956699871
10  15.048692390077921
20  105.98527934204294
40  159.67329917703637
100 258.8574608904579
};
\addlegendentry{EnKF $10^6$}

% \addplot [semithick, darkorange25512714, mark=*, mark size=1, mark options={solid}]
\addplot[method PF 1e6, with marks every]
table {%
4   100.75090938916307
10  364.3423469229075
20  539.915313953252
40  941.9727584683784
100 1622.8431585704996
};
\addlegendentry{PF $10^6$}

% --- reference slope guides on log–log ---
\pgfmathsetmacro{\dref}{10}
\pgfmathsetmacro{\yref}{20}

% O(d^{1/2})
\addplot [semithick, black, dashed, domain=4:100, samples=2]
  {\yref * (x/\dref)^(0.5)};
\addlegendentry{$O(d^{1/2})$}

% O(d)
\addplot [semithick, black, mark=*, mark size=1, mark options={solid}]
table{
0.01   100.75090938916307
};
\addlegendentry{$O(d)$}

\nextgroupplot[
  legend cell align={left},
  legend style={
    fill opacity=0.8,
    draw opacity=1,
    text opacity=1,
    at={(0.03,0.97)},
    anchor=north west,
    draw=none,
    nodes={scale=1.0, transform shape},
    legend columns=-1
  },
  tick align=outside,
  tick pos=left,
  title={NLL},
  x grid style={darkgray176},
  xlabel={$d$},
  xmajorgrids,
  xmin=3, xmax=110,
  xmode=log,
  log basis x=2,
  log ticks with fixed point,
  xtick={4,10,20,40,100},
  xtick style={color=black},
  y grid style={darkgray176},
  ymajorgrids,
  ymode=log,
  log basis y=10,
  log ticks with fixed point,
  ymin=1e1, ymax=6e2,
  ytick style={color=black},
  legend to name={legend: Lorenz highdim NLL group},
]

% --- your data series (unchanged) ---
% \addplot [semithick, steelblue31119180, mark=*, mark size=1, mark options={solid}]
\addplot[method LogBSDEF, with marks every]
table {%
4   13.29342128632654
10  39.59141707609639
20  63.06306888509524
40  94.23963266917845
100 125.4939480161006
};
\addlegendentry{LogBSDEF}

% \addplot [semithick, crimson2143940, mark=*, mark size=1, mark options={solid}]
\addplot[method EnKF 1e6, with marks every]
table {%
4   22.225008648561484
10  36.92781456025328
20  66.42809670918811
40  134.45547946143455
100 331.7637193923017
};
\addlegendentry{EnKF $10^6$}

% \addplot [semithick, darkorange25512714, mark=*, mark size=1, mark options={solid}]
\addplot[method PF 1e6, with marks every]
table {%
4   238.91105293017375
10  268.4033890513059
20  297.7422216686388
40  358.9642104332838
100 423.94259154363357
};
\addlegendentry{PF $10^6$}

% --- reference slope guide on log–log ---
\pgfmathsetmacro{\dref}{10}
\pgfmathsetmacro{\yref}{50}

% O(d)
\addplot [semithick, black, domain=4:100, samples=2]
  {\yref * (x/\dref)};
\addlegendentry{$O(d)$}

\end{groupplot}

\end{tikzpicture}
    \begin{tikzpicture}[overlay]
    \node at (0,0) {$ $};
    \node at (-6.0,-0.3) {\ref*{legend: Lorenz highdim MAE group}};
    \end{tikzpicture}
    \vspace{1em}
    \caption{The averaged MAE and NLL metrics evaluated over state dimension $d = [4,10,20,40,100]$. At each time step the NLL value is capped at $|\log(10^{-200})| \approx 460$.}
    \vspace{-0.5em}
    \label{fig: Lorenz highdim}
\end{figure}

\subsection{Computational efficiency}
\label{section: computional efficiency}
Finally, we study the inference time of the methods and how it scales with the state dimension. We display the time required to i) estimate one whole state trajectory by the first moment for each observation time and ii) evaluate each filtering density in 1000 points. The deep learning-based filters were tested using a single NVIDIA A40 GPU with 48 GB of VRAM. The PF and EnKF were run using an Intel(R) Xeon(R) Platinum 8358 CPU with 32 cores.

We remark that the runtime for all the deep density methods primarily depends on the dimension of the problem (which more or less determines the size of the neural network) as well as the sampling method and the number of samples used to compute the first moment or normalizing constant. Hence, we show the runtime when fixing $T = 1$, $K=10$ and $d=d'$ for the Ornstein--Uhlenbeck case, and argue that the other problems have more or less the same computational cost. In this setting we show the inference times in Figure \ref{fig: runtime}. As can be seen, BSDEF (whose runtime is equal to that of DSF, LogBSDEF, and LogDSF) scales well with the state dimension, with runtimes remaining rather constant regardless of the sampling method (abbreviated by I-EKF for the EKF importance distribution and I-G for the observation independent Gaussian importance distribution, see Supplementary~Materials~\ref{SM: Normalization} for details). In comparison, the benchmark filters explode in estimation time with higher $d$ and their density calculation time is consistently many magnitudes higher due to the expensive kernel density estimate that we utilize to obtain densities (see Supplementary~Materials \ref{SM: Evaluation} for details). 

In Figure~\ref{fig: training_time} we report the total time to perform training, state estimation, and density evaluation (with $1000$ points per observation time) on two examples, as a function of the number of observation sequences. We include the best-performing deep density based methods alongside a particle filter and an ensemble Kalman filter. For the $100$-dimensional Lorenz-96 case, the classical methods do not yield accurate estimates. Their curves are shown only to illustrate computational cost. The trade-off is clear: deep density based methods incur a one time training cost but offer low per-sample inference time, whereas particle-based methods avoid training but have higher per-sample cost. In terms of computation time, the deep density based methods become advantageous once more than about $1300$ (bistable) or $430$ (Lorenz-96) sequences are processed.

\begin{figure}[h]
    \centering
    % This file was created with tikzplotlib v0.10.1.
\begin{tikzpicture}

\definecolor{crimson2143940}{RGB}{214,39,40}
\definecolor{darkgray176}{RGB}{176,176,176}
\definecolor{darkorange25512714}{RGB}{255,127,14}
\definecolor{darkturquoise23190207}{RGB}{23,190,207}
\definecolor{forestgreen4416044}{RGB}{44,160,44}
\definecolor{goldenrod18818934}{RGB}{188,189,34}
\definecolor{gray127}{RGB}{127,127,127}
\definecolor{lightgray204}{RGB}{204,204,204}
\definecolor{mediumpurple148103189}{RGB}{148,103,189}
\definecolor{orchid227119194}{RGB}{227,119,194}
\definecolor{sienna1408675}{RGB}{140,86,75}
\definecolor{steelblue31119180}{RGB}{31,119,180}

\begin{groupplot}[group style={group size=2 by 1,horizontal sep=1.5cm,
    vertical sep=1.5cm}, width=7.0cm,
  height=5.0cm]
\nextgroupplot[
legend cell align={left},
legend style={
  fill opacity=0.8,
  draw opacity=1,
  text opacity=1,
  at={(0.03,0.97)},
  anchor=north west,
  draw=none,
  nodes={scale=1.0, transform shape},
  legend columns=2,
  transpose legend
},
tick align=outside,
tick pos=left,
title={Estimation Time},
x grid style={darkgray176},
xlabel={$d$},
xmajorgrids,
xmin=0.445869142790845, xmax=100.15454867838,
xminorgrids,
xtick style={color=black},
xmin=0.9, xmax=110,         % nicer margins on a log scale
xmode=log,
log basis x=2,
log ticks with fixed point,
xtick={1,2,4,10,20,40,100},
y grid style={darkgray176},
ytick={0.01,0.1,1,10,100},
ylabel={Time (s)},
ymajorgrids,
ymode=log,
log basis y=2,
ymin=-0.0162472780164953, ymax=300.706815206030834321,
yminorgrids,
ytick style={color=black},
legend to name={legend: runtime},
]

% BSDEF I–EKF  (use method BSDEF)
% \addplot [semithick, crimson2143940, mark=*,  mark size=1, mark options={solid}]
\addplot[method BSDEF, with marks]
table {%
1 0.0822582836151123
2 0.08610183358192444
4 0.0887244467735291
10 0.0978767921924591
20 0.102883826255798
40 0.10544020676612854
100 0.148568474769592
};
\addlegendentry{BSDEF I-EKF}

% BSDEF I–G  (use method LogBSDEF)
% \addplot [semithick, sienna1408675, mark=*,  mark size=1, mark options={solid}]
\addplot[method LogBSDEF, with marks]
table {%
1 0.0188338766098022
2 0.018940519332885743
4 0.0222576761245728
10 0.0304657287597656
20 0.0309862058162689
40 0.03343233752250672
100 0.0712657968997955
};
\addlegendentry{BSDEF I-G}

% ENKF 1e4
% \addplot [semithick, steelblue31119180, mark=*,  mark size=1, mark options={solid}]
\addplot[method EnKF 1e4, with marks]
table {%
1 0.00807343935966492
2 0.016027955532073975
4 0.0202659697532654
10 0.0243203763961792
20 0.407733315944672
40 0.48563116216659546
100 1.20042042303085
};
\addlegendentry{EnKF $10^4$}

% ENKF 1e5
% \addplot [semithick, darkorange25512714, mark=*,  mark size=1, mark options={solid}]
\addplot[method EnKF 1e5, with marks]
table {%
1 0.0224012362957001
2 0.038641357898712154
4 0.0632707281112671
10 0.162620323181152
20 3.89580481863022
40 4.727338916301727
100 11.8950538492203
};
\addlegendentry{EnKF $10^5$}

% ENKF 1e6
% \addplot [semithick, forestgreen4416044, mark=*,  mark size=1, mark options={solid}]
\addplot[method EnKF 1e6, with marks]
table {%
1 0.09788
2 0.25968
4 0.59657
10 1.83865
20 40.401354
40 47.9983
100 112.4534
};
\addlegendentry{EnKF $10^6$}

% PF 1e4
% \addplot [semithick, cyan, mark=*,  mark size=1, mark options={solid}]
\addplot[method PF 1e4, with marks]
table {%
1 0.0192896375656128
2 0.021834341526031493
4 0.0229717643260956
10 0.0259979720115662
20 0.231718129634857
40 0.27015773153305056
100 0.678173982620239
};
\addlegendentry{PF $10^4$}

% PF 1e5
% \addplot [semithick, yellow, mark=*,  mark size=1, mark options={solid}]
\addplot[method PF 1e5, with marks]
table {%
1 0.163513789892197
2 0.15812388610839845
4 0.165065482378006
10 0.201421261787415
20 2.15993277168274
40 2.6336552429199216
100 6.59376759195328
};
\addlegendentry{PF $10^5$}

% PF 1e6
% \addplot [semithick, mediumpurple148103189, mark=*,  mark size=1, mark options={solid}]
\addplot[method PF 1e6, with marks]
table {%
1 2.7967
2 2.7854
4 2.96644
10 3.51111
20 23.92237
40 28.1656
100 66.0933
};
\addlegendentry{PF $10^6$}

\nextgroupplot[
legend cell align={left},
legend style={fill opacity=0.8, draw opacity=1, text opacity=1, draw=lightgray204},
tick align=outside,
tick pos=left,
title={Density Calc. Time},
x grid style={darkgray176},
xlabel={$d$},
xmajorgrids,
xmin=0.445869142790845, xmax=100.15454867838,
xminorgrids,
xtick style={color=black},
xmin=0.9, xmax=110,         % nicer margins on a log scale
xmode=log,
log basis x=2,
log ticks with fixed point,
xtick={1,2,4,10,20,40,100},
y grid style={darkgray176},
ymajorgrids,
ymode=log,
log basis y=2,
ymin=-0.0162472780164953, ymax=2000.706815206030834321,
ytick={0.01,0.1,1,10,100,1000},
yminorgrids,
ytick style={color=black},
]

% BSDEF Importance 1e4  (map to method BSDEF)
% \addplot [semithick, crimson2143940, mark=*,  mark size=1, mark options={solid}]
\addplot[method BSDEF, with marks]
table {%
1 0.0840211224555969
2 0.0895618224143982
4 0.0884912419319153
10 0.0955559587478638
20 0.101902439594269
40 0.105348887443542
100 0.142997376918793
};
%\addlegendentry{DSF/BSDEF Importance $10^4$}

% LogBSDEF IS 1e4  (map to method LogBSDEF)
% \addplot [semithick, sienna1408675, mark=*,  mark size=1, mark options={solid}]
\addplot[method LogBSDEF, with marks]
table {%
1 0.0210011672973633
2 0.0217149925231934
4 0.0217698264122009
10 0.0287273693084717
20 0.0294194102287292
40 0.0318874049186707
100 0.0679280352592468
};
%\addlegendentry{DSF/BSDEF IS $10^4$}

% ENKF 1e4
% \addplot [semithick, steelblue31119180, mark=*,  mark size=1, mark options={solid}]
\addplot[method EnKF 1e4, with marks]
table {%
1 1.38944756031036
2 1.14627522945404
4 1.23352312088013
10 1.67014675617218
20 2.81462216615677
40 4.97886157989502
100 12.4661632204056
};
%\addlegendentry{EnKF $10^4$}

% ENKF 1e5
% \addplot [semithick, darkorange25512714, mark=*,  mark size=1, mark options={solid}]
\addplot[method EnKF 1e5, with marks]
table {%
1 16.9234394931793
2 16.4461768937111
4 14.4525819635391
10 16.5917891144752
20 28.0417708706856
40 49.513944747448
100 126.097231492996
};
%\addlegendentry{EnKF $10^5$}

% ENKF 1e6
% \addplot [semithick, forestgreen4416044, mark=*,  mark size=1, mark options={solid}]
\addplot[method EnKF 1e6, with marks]
table {%
1 174.9903
2 230.9410
4 221.7002
10 168.6244
20 280.5036
40 491.0797
100 1290.08
};
%\addlegendentry{EnKF $10^6$}

% PF 1e4
% \addplot [semithick, cyan, mark=*,  mark size=1, mark options={solid}]
\addplot[method PF 1e4, with marks]
table {%
1 1.38898206472397
2 1.12093079090118
4 1.22891935825348
10 1.66692132949829
20 2.64525503635406
40 5.57793086051941
100 10.8449607872963
};
%\addlegendentry{PF $10^4$}

% PF 1e5
% \addplot [semithick, yellow, mark=*,  mark size=1, mark options={solid}]
\addplot[method PF 1e5, with marks]
table {%
1 16.9555643391609
2 15.383254032135
4 13.3015938329697
10 16.5907232880592
20 26.2638012814522
40 50.7173513150215
100 108.784478769302
};
%\addlegendentry{PF $10^5$}

% PF 1e6
% \addplot [semithick, mediumpurple148103189, mark=*,  mark size=1, mark options={solid}]
\addplot[method PF 1e6, with marks]
table {%
1 177.5215
2 224.1256
4 190.3013
10 165.1656
20 263.2891
40 485.5005
100 1146.39
};
%\addlegendentry{PF $10^6$}

\end{groupplot}

%\node at (0,0) {$ $};
%\node at (5,-1.6) {\ref*{legend: runtime}};

\end{tikzpicture}
    \begin{tikzpicture}[overlay]
    \node at (0,0) {$ $};
    \node at (-6,-0.6) {\ref*{legend: runtime}};
    \end{tikzpicture}
    \vspace{2.5em}
    \caption{On the left, we display the average time for estimating one whole trajectory in the Ornstein--Uhlenbeck case. On the right, we display the average time for evaluating filtering densities for all observation times in $1000$ points. Note that this includes the time to obtain normalization constants for the BSDEF. In both plots, the time it takes to propagate particles for the EnKF and PF is included. For BSDEF, we fix the number of samples to $10^4$ used to compute the first moment and normalizing constant, although it is normally possible to use fewer samples without sacrificing much performance.}
    \vspace{-0.5em}
    \label{fig: runtime}
\end{figure}
\begin{figure}[h!]
    \centering
    % This file was created with tikzplotlib v0.10.1.
\begin{tikzpicture}

\definecolor{crimson2143940}{RGB}{214,39,40}
\definecolor{darkgray176}{RGB}{176,176,176}
\definecolor{darkorange25512714}{RGB}{255,127,14}
\definecolor{forestgreen4416044}{RGB}{44,160,44}
\definecolor{gray127}{RGB}{127,127,127}
\definecolor{lightgray204}{RGB}{204,204,204}
\definecolor{mediumpurple148103189}{RGB}{148,103,189}
\definecolor{orchid227119194}{RGB}{227,119,194}
\definecolor{sienna1408675}{RGB}{140,86,75}
\definecolor{steelblue31119180}{RGB}{31,119,180}

\begin{groupplot}[group style={group size=2 by 1,horizontal sep=1.5cm,
    vertical sep=1.5cm}, width=7.0cm,
  height=5.0cm]

%% Bistable
\nextgroupplot[
%width=0.35*6.028in,
%height=0.38*4.754in,
legend cell align={left},
legend style={
  fill opacity=0.8,
  draw opacity=1,
  text opacity=1,
  at={(0.03,0.97)},
  anchor=north west,
  draw=none,
  nodes={scale=1.0, transform shape},
  legend columns=-1
},
tick align=outside,
tick pos=left,
title={Bistable},
x grid style={darkgray176},
xlabel={Processed observations},
xmajorgrids,
xmin=0.55, xmax=200000.045,
xtick style={color=black},
xmode=log,
log basis x=10, 
ymode=log,
log basis y=10, 
y grid style={darkgray176},
ylabel={Time (s)},
ymajorgrids,
ymin=1, ymax=30000000,
ytick style={color=black},
legend to name={legend: training to inference},
xtick={1,10,100,1000,10000,100000},
xticklabels={1,10,$10^{2}$,$10^{3}$,$10^{4}$,$10^{5}$},
]
\addplot [method LogBSDEF]
table {%
0.00001 0.00001
1 219600
10 219600
100 219602
1000 219621
100000 221700
1000000 221700
};
\addlegendentry{LogBSDEF}

\addplot [method LogDSF, dashed]
table {%
0.00001 0.00001
1 230400
100 230402
10000 230610
1000000 232500
};
\addlegendentry{LogDSF}

\addplot [method EnKF 1e6]
table {%
0.001 0.175
%1 175
10 1750
1000 175000
10000  1750000
100000  17500000
1000000 175000000
};
\addlegendentry{EnKF $10^6$}

\addplot [method PF 1e6, dashed]
table {%
0.001 0.177
1 177
100 17700
10000 1770000
1000000 177000000
};
\addlegendentry{PF $10^6$}

%% L96 100D
\nextgroupplot[
%width=0.35*6.028in,
%height=0.38*4.754in,
legend cell align={left},
legend style={
  fill opacity=0.8,
  draw opacity=1,
  text opacity=1,
  at={(0.03,0.97)},
  anchor=north west,
  draw=none,
  nodes={scale=1.0, transform shape},
  legend columns=-1
},
tick align=outside,
tick pos=left,
title={Lorenz-96 (100D)},
x grid style={darkgray176},
xlabel={Processed observations},
xmajorgrids,
xmin=0.55, xmax=200000.045,
xtick style={color=black},
xmode=log,
log basis x=10, 
ymode=log,
log basis y=10, 
y grid style={darkgray176},
ymajorgrids,
ymin=1, ymax=30000000,
ytick style={color=black},
legend to name={legend that is not used},
xtick={1,10,100,1000,10000,100000},
xticklabels={1,10,$10^{2}$,$10^{3}$,$10^{4}$,$10^{5}$},
]
\addplot [method LogBSDEF]
table {%
0.00001 0.00001
1 496800
10 496801
10 496800 + 0.0679
100 496807
1000 496868
10000 497480
100000 503600
1000000 564800
};
\addlegendentry{LogBSDEF}

\addplot [method EnKF 1e6]
table {%
0.001 1.290
10 12900
1000   1290000
100000 129000000
};
\addlegendentry{EnKF $10^6$}

\addplot [method PF 1e6, dashed]
table {%
0.001 1.146
1 1146
100    114600
10000  11460000
100000 114600000
};
\addlegendentry{PF $10^6$}

\end{groupplot}

\end{tikzpicture}
    \begin{tikzpicture}[overlay]
    \node at (0,0) {$ $};
    \node at (-6.7,-0.5) {\ref*{legend: training to inference}};
    \end{tikzpicture}
    \vspace{2em}
    \caption{The computational time, from initializing each method to evaluation of $1000$ spatial points in each observation time, over increasing number of sequences. The intersection between the PFs and deep density methods' computational time occurs at $1300$ samples for the bistable example and at $430$ samples for the $100$-dimensional Lorenz-96 example.}
    \vspace{-0.5em}
    \label{fig: training_time}
\end{figure}

\section{Conclusion and discussion}
\label{section: conclusion}
This work benchmarked deep filtering methods, the Deep Splitting Filter (DSF) and the deep Backward Stochastic Differential Equation Filter (BSDEF), together with their logarithmic variants for nonlinear state estimation in high-dimensional systems. By formulating prediction and update probabilistically, the methods can be trained and allow efficient evaluation.

Across all testbeds, a clear pattern emerges. On linear problems with reliable ground truth, the methods recover high accuracy and remain stable in long horizons and high dimensions. The logarithmic formulations improve numerical robustness and consistently reduce error, with LogBSDEF giving the strongest results. In addition, the log-scale formulation yields positivity-preserving density approximations, eliminating inconsistent negative density estimates and contributing to stable training and evaluation in high dimensions.
On strongly nonlinear dynamics, including the Schl\"ogl model and stochastic Lorenz-96 system, the deep density filters remain numerically stable and often competitive, whereas particle and ensemble methods tend to degrade or become costly when the underlying dimension increases. The LogBSDEF exhibits only mild growth in inference time with dimension, supporting its use in high-dimensional settings, unlike particle filters and the ensemble Kalman filter, which inherently incur the curse of dimensionality. 
The difference in computational time is about two to five orders of magnitude.

Several limitations remain. Performance depends on network architecture and sampling choices, and the selection of the initialization density $q_0$ can influence stability in stiff regimes. Our experiments use Gaussian observation noise and fixed time discretizations; extending to alternative likelihoods and adaptive stepping is left for future work. Finally, while the methods can achieve high accuracy and stability across many settings, they require careful hyperparameter tuning which remains a practical obstacle. Supplementary~Materials~\ref{SM: Training} documents the configurations, which confirms the variety of parameters used. 

At the same time, the experiments with the LSTM-based model in the long-horizon setting indicate that more advanced neural network architectures can help improve performance. However, we tried a few other architectures that did not make it into the paper, among them transformer-based and other LSTM architectures, which, for these examples and configurations, did not improve performance. Further work in this direction is needed in order to find optimal architectures.

\subsection*{Acknowledgements}
The authors would like to thank Adam Andersson, Stig Larsson, and Ruben Seyer for valuable input and discussions on this manuscript.
The work of {K.B.} was supported by the Wallenberg AI, Autonomous Systems and Software Program (WASP) funded by the Knut and Alice Wallenberg Foundation. The work of {F.R.} was funded by the Swedish Electromobility Centre (SEC) and partially supported by WASP. The computations were enabled by resources provided by the National Academic Infrastructure for Supercomputing in Sweden (NAISS) at Chalmers e-Commons partially funded by the Swedish Research Council through grant agreement no. 2022-06725.

\appendix
\section{Proof of Theorem~\ref{theorem: log equation}}
\label{appendix: Proofs}
For simplicity, we hide $t$, $x$, $o$, and $k$ in the notation and write $p=p_k(t,x,o_{1:k})$.
The initial condition follows by insertion and simplification. It remains to derive the log-transformed equation from the original Fokker--Planck equation. 
By taking the time derivative in $v=-\log p$, we obtain
\begin{align}
    \label{proof eq: diff pk vk relation 0}
    \frac{\partial v}{\partial t}
    &=
    -\frac{1}{p}
    \frac{\partial p}{\partial t}
    ,
\end{align}
or equivalently 
\begin{align} \label{proof eq: diff pk vk relation}
    \frac{\partial p}{\partial t}
    &=
    -
    p
    \frac{\partial v}{\partial t}
    .
\end{align}
Analogous relations hold for the spatial derivatives. We consider $Ap$ and substitute $\frac{\partial p}{\partial x_i}$ by \eqref{proof eq: diff pk vk relation} to get
\begin{align*}
    Ap
    &= 
    \frac{1}{2}\sum_{i,j=1}^d a_{ij}\,
    \frac{\partial^2 p}{\partial x_i \partial x_j} 
    + 
    \sum_{i=1}^d 
    \mu_i \, 
    \frac{\partial p}{\partial x_i}
    =
    -
    \frac{1}{2}\sum_{i,j=1}^d a_{ij}\,
    \frac{\partial}{\partial x_i } 
    \bigg(
    p
    \frac{\partial v}{\partial x_j} 
    \bigg)
    -
    \sum_{i=1}^d 
    \mu_i \, 
    p
    \frac{\partial v}{\partial x_i}.
\end{align*}
Next, by the product rule and \eqref{proof eq: diff pk vk relation}, we obtain 
\begin{align}
\begin{split}
\label{proof eq: Ap Av relation}
    Ap
    &=
    \frac{1}{2}\sum_{i,j=1}^d a_{ij}\, 
    p
    \bigg(
    \frac{\partial v}{\partial x_i }
    \frac{\partial v}{\partial x_j} 
    \bigg)
    -
    \frac{1}{2}\sum_{i,j=1}^d a_{ij}\, 
    p
    \bigg(
    \frac{\partial^2 v}{\partial x_i \partial x_j} 
    \bigg)
    -
    \sum_{i=1}^d 
    \mu_i \, 
    p
    \frac{\partial v}{\partial x_i}
    \\
    &=
    p
    \Big(
    \frac{1}{2}
    \|
    \sigma^\top
    \nabla v
    \|^2
    -
    Av
    \Big).    
\end{split}
\end{align}
By inserting \eqref{proof eq: Ap Av relation} into the differentiated form of the Fokker--Planck equation, we get
\begin{align} 
\begin{split}
    \label{proof eq: diff FP-prediction}
    \frac{\partial }{\partial t}
    p
    &=
    A
    p
    +
    f(
    p,
    \nabla
    p
    )
    =
    p
    \Big(
    \frac{1}{2}
    \|
    \sigma^\top
    \nabla v
    \|^2
    -
    Av
    \Big)
    +
    f(
    p,
    \nabla
    p
    )
    . 
\end{split}
\end{align}
We continue by inserting \eqref{proof eq: diff FP-prediction} into \eqref{proof eq: diff pk vk relation 0} to obtain
\begin{align*}
    \frac{\partial}{\partial t}
    v
    &=
    -\frac{1}{p}
    \Big(
    p
    \Big(
    \frac{1}{2}
    \|
    \sigma^\top
    \nabla v
    \|^2
    -
    Av
    \Big)
    +
    f(
    p,
    \nabla
    p
    )
    \Big)
    \\
    &=
    Av
    -
    \frac{1}{2}
    \|
    \sigma^\top
    \nabla v
    \|^2
    -
    \frac{1}{p}
    f(
    p,
    \nabla
    p
    ).
\end{align*}
Since $f$ is linear, we obtain
\begin{align*}
    \frac{\partial}{\partial t}
    v
    &=
    Av
    -
    \frac{1}{2}
    \|
    \sigma^\top
    \nabla v
    \|^2
    -
    f(
    1,
    -
    \nabla
    v
    ).
\end{align*}
The proof is finished by recalling the definition of $f_{\text{log}}$ from the theorem statement. 

%\bibliographystyle{siamplain}
%\bibliography{references}

%\putbib
%\end{bibunit}
\printbibliography[title={References}]
\end{refsection}

\clearpage
\printaddresseshere
\makeatletter
\let\enddoc@text\relax
\makeatother
\clearpage

\setcounter{page}{1}

\markboth{SUPPLEMENTARY MATERIALS: FILTERING WITH DEEP DENSITY METHODS}{K. B\r{A}GMARK AND F. RYDIN}

\thispagestyle{plain}

%\clearpage
\begin{center}
    {\large \textbf{SUPPLEMENTARY MATERIALS: HIGH-DIMENSIONAL BAYESIAN FILTERING THROUGH DEEP DENSITY APPROXIMATION}}\\[1.5em]
    {\small KASPER BÅGMARK AND FILIP RYDIN}
\end{center}
\vspace{2em}

% reset numbering for supplement
\setcounter{section}{0}
\setcounter{subsection}{0}
\setcounter{subsubsection}{0}
\setcounter{equation}{0}
\setcounter{figure}{0}
\setcounter{table}{0}
\setcounter{theorem}{0}
\setcounter{lemma}{0}
\setcounter{proposition}{0}
\setcounter{corollary}{0}
\setcounter{definition}{0}
\setcounter{remark}{0}
\setcounter{assumption}{0}
\setcounter{setting}{0}
\setcounter{claim}{0}

\renewcommand{\thesection}{SM\arabic{section}}
\renewcommand{\thesubsection}{\thesection.\arabic{subsection}}
\renewcommand{\thesubsubsection}{\thesubsection.\arabic{subsubsection}}
\renewcommand{\theequation}{SM\arabic{section}.\arabic{equation}}
\renewcommand{\thefigure}{SM\arabic{figure}}
\renewcommand{\thetable}{SM\arabic{table}}

%\begin{bibunit}
\begin{refsection}

\section{Classical methods}
\label{SM: Classical methods}
In this section we recall the classical methods that we ran our comparisons against. This includes a bootstrap particle filter, the extended Kalman filter, and an ensemble Kalman filter. 
\subsection{Particle filter}
\label{section: particle filter}
The classical bootstrap particle filter, a type of sequential Monte Carlo method, was introduced in \cite{gordon1993novel}, and provides an approximation of the nonlinear filter by propagating and reweighting an ensemble of particles. The method is based on $M\in \mathbb{N}$ weighted particles $(\widetilde{X}_k^{(i)},V_k^{(i)})_{i=1}^M$, with $\widetilde{X}_k^{(i)}$ denoting the particle and $V_k^{(i)}$ its weight, approximating the posterior distribution at each observation time $t_k$, $k=1,\dots,K$.

The filter is initialized by sampling each $\widetilde{X}_0^{(i)}$, $i=1,\dots,M$, according $\pi_0$, and setting $V_0^{(i)} = \frac{1}{M}$ for $i=1,\dots,M$. Next, given $(\widetilde{X}_{k-1}^{(i)},V_{k-1}^{(i)})_{i=1}^M$ that approximates $p_{k-1}(t_{k-1})$ with an empirical distribution, the bootstrap particle filter proceeds by:
\begin{enumerate}
    \item \textbf{Prediction:} Propagate each particle forward according to $F$ approximating the dynamics of the signal process $S$, from \eqref{eq: state}, by
    \begin{align*}
        \widehat{X}_k^{(i)} 
        &=
        F
        \bigl(
            \widetilde{X}_{k-1}^{(i)},\,
            \Delta W_k^{(i)}
        \bigr)
        ,\quad
        i=1,\dots,M
        .
    \end{align*}
    Here $\Delta W_k^{(i)}$ is the Brownian motion increment for the $i$'th particle at time $t_k$. For instance we can let $F$ be the forward operator given by the Euler--Maruyama discretization defined from \eqref{eq: Euler--Maruyama SDE}.
    \item \textbf{Weighting:} Compute importance weights proportional to the likelihood of the observation $O_k$ given the predicted particles:
    \begin{align*}
        \widetilde{V}_k^{(i)} = 
        V_{k-1}^{(i)} \, L(O_k,\widehat{X}_k^{(i)}),
        \quad
        i=1,\dots,M.
    \end{align*}
    Normalize weights:
    \begin{align*}
        \bar{V}_k^{(i)}
        =
        \frac{\widetilde{V}_k^{(i)}}
        {\sum_{j=1}^M \widetilde{V}_k^{(j)}}
        ,\quad
        i=1,\dots,M.
    \end{align*}
    \item \textbf{Resampling:} Generate a new set of particles $(\widetilde{X}_k^{(i)},V_k^{(i)})_{i=1}^M$ by sampling with replacement from $\{\widehat{X}_k^{(i)}\}_{i=1}^M$ according to the new normalized weights $\{\bar{V}_k^{(i)}\}_{i=1}^M$, and reset weights to $V_k^{(i)} = 1/M$.
\end{enumerate}
Under standard regularity assumptions (bounded likelihood, nondegenerate resampling), the empirical measure 
\begin{align*}
    \overline{p}_k^M 
    &= 
    \frac{1}{M}
    \sum_{i=1}^M 
    \delta_{\widetilde{X}_k^{(i)}}
\end{align*}
is a consistent estimator of the filtering distribution in the sense that, for any bounded test function $\varphi$,
\begin{align*}
    \overline{p}_k^M(\varphi)
    &=
    \frac{1}{M}
    \sum_{i=1}^M
    \varphi
    (
    \widetilde{X}_k^{i}
    )
    \xrightarrow[M\to\infty]{a.s.}
    \mathbb{E}
    \left[
    \varphi(S_{t_k}) 
    \mid 
    O_{1:k}
    \right].
\end{align*}
Moreover, a central limit theorem holds with $\sqrt{M}$-rate convergence and an explicit asymptotic variance, see \cite{chopin2004central, crisan2002survey, moral2004feynman}. The bootstrap particle filter thus provides a sequential Monte Carlo approximation of the nonlinear filter that converges almost surely as $M\to\infty$, but suffers from the curse of dimensionality, as its variance grows exponentially with the dimension of the state space \cite{bickel2008sharp}.

\subsection{Extended Kalman filter}
For completeness we briefly recall the main steps of the continuous-discrete extended Kalman filter, which provides a Gaussian approximation of the filtering distribution. The EKF propagates the conditional mean and covariance through a local linearization of the dynamics and a Kalman-like update at observation times. We let $(m_k,P_k)_{k=0}^K$ denote the mean and covariance parameterizing the EKF at each time $t_k$. If $\pi_0$ is a Gaussian density, we simply initiate $m_0$ and $P_0$ with the true parameters from $\pi_0$. The method recursively approximates the prediction step \eqref{eq: FP-prediction} and the update step \eqref{eq: FP-update}.
\begin{enumerate}
    \item \textbf{Prediction:}
        Given the mean and covariance $(m_{k-1},P_{k-1})$ at $t_{k-1}$, we integrate the moment equations, initialized at $\widehat{m}(t_{k-1}) = m_{k-1}$ and $\widehat{P}(t_{k-1}) = P_{k-1}$, over $t\in[t_{k-1},t_k]$ 
        \begin{align*}
            \frac{\dd \widehat{m}(t)}{\dd t}
            &= 
            \mu(\widehat{m}(t)),
            &
            \frac{\dd \widehat{P}(t)}{\dd t}
            &=
            A(t) \widehat{P}(t) 
            + 
            \widehat{P}(t) A(t)^\top 
            + 
            \sigma
            \big(\widehat{m}(t)
            \big)
            \sigma
            \big(\widehat{m}(t)
            \big)^\top
            ,
        \end{align*}
        where $A(t) \defeq \mathrm{D} \mu(\widehat{m}(t))$ is the Jacobian of the drift $\mu$ evaluated at $\widehat{m}(t)$. We denote the moment predictions at the terminal time by 
        \begin{align*}
            m_{k|k-1} = \widehat{m}(t_k), 
            \qquad
            P_{k|k-1} = \widehat{P}(t_k).
        \end{align*}
    \item \textbf{Update:}
        Let $H_k \defeq \mathrm{D} h(m_{k|k-1})$ denote the Jacobian of the observation function $h$ evaluated at the predicted mean. Then the Kalman gain is
        \begin{align*}
            K_k &= 
            P_{k|k-1} H_k^\top 
            \big(H_k P_{k|k-1} H_k^\top + R\big)^{-1},
        \end{align*}
        where we recall that $R$ is the observation noise covariance. The posterior mean and covariance become
        \begin{align*}
            m_k &= 
            m_{k|k-1} + 
            K_k \big(O_k - h(m_{k|k-1})\big),
            &
            P_k &= 
            (I - K_k H_k) P_{k|k-1}.
        \end{align*}
\end{enumerate}
The EKF thus produces a Gaussian $\mathcal{N}(x \mid m_k,P_k)$ approximating the filter $p_k(t_k,x)$. This distribution can be used as a closed-form approximation of the filter, but also as a proposal distribution for Monte Carlo methods such as importance sampling. This is utilized in the high-dimensional examples and further explained in Supplementary~Materials~\ref{SM: Normalization}.

\subsection{Ensemble Kalman filter}
The ensemble Kalman filter \cite{evensen1994sequential} is a Monte Carlo implementation of the Kalman filter methodology, designed to make Kalman-type filtering feasible for high-dimensional, possibly nonlinear models. Rather than propagating a mean and covariance matrix, the EnKF maintains an ensemble of $M$ particles $\{\widetilde{X}_{k-1}^{(i)}\}_{i=1}^M$ that approximates the filtering distribution at time $t_{k-1}$. The method alternates between a prediction step, where each ensemble member is advanced forward in time, and an update step, where information from the new observation $O_k$ is assimilated. We present here the stochastic or perturbed-observations variant of the EnKF, which preserves the posterior covariance structure in expectation \cite{burgers1998analysis,katzfuss2016understanding}.
\begin{enumerate}
    \item \textbf{Prediction:}
    Each ensemble member is propagated according to some forward operator $F$ that approximates the dynamics of $S$, for example using the Euler--Maruyama scheme \eqref{eq: Euler--Maruyama SDE}, to form the forecast ensemble
    \begin{align*}
        \widehat{X}_k^{(i)}
        &=
        F
        \bigl(
            \widetilde{X}_{k-1}^{(i)},\,
            \Delta W_k^{(i)}
        \bigr),
        \qquad i=1,\dots,M.
    \end{align*}
    Here $\Delta W_k^{(i)}$ is the Brownian motion increment for the $i$'th ensemble member at time $t_k$. The forecast sample mean and covariance are then computed as
    \begin{align*}
        \overline{x}_k
        &=
        \frac{1}{M}
        \sum_{i=1}^M
        \widehat{X}_k^{(i)},
        &
        P^x_k
        &=
        \frac{1}{M-1}
        \sum_{i=1}^M
        \bigl(
            \widehat{X}_k^{(i)}
            - 
            \overline{x}_k
        \bigr)
        \bigl(
            \widehat{X}_k^{(i)}
            - 
            \overline{x}_k
        \bigr)^\top.
    \end{align*}
    
    \item \textbf{Update:}
    Draw independent perturbations $\varepsilon_k^{(i)}\sim \mathcal{N}(0,R)$ and define the predicted observations, and their sample mean, by
    \begin{align*}
        \widehat{Y}_k^{(i)}
        &=
        h(\widehat{X}_k^{(i)})
        +
        \varepsilon_k^{(i)}
        ,
        \quad
        \overline{y}_k
        =
        \frac{1}{M}
        \sum_{i=1}^M
        \widehat{Y}_k^{(i)}
        .
    \end{align*}
    From the forecast ensemble $(\widehat{X}_k^{(i)},\widehat{Y}_k^{(i)})_{i=1}^M$ we have the covariance matrices
    \begin{align*}
        P^{xy}_k
        &=
        \frac{1}{M-1}
        \sum_{i=1}^M
        \bigl(
            \widehat{X}_k^{(i)} - \overline{x}_k
        \bigr)
        \bigl(
            \widehat{Y}_k^{(i)} - \overline{y}_k
        \bigr)^\top,
        \\
        P^{yy}_k
        &=
        \frac{1}{M-1}
        \sum_{i=1}^M
        \bigl(
            \widehat{Y}_k^{(i)} - \overline{y}_k
        \bigr)
        \bigl(
            \widehat{Y}_k^{(i)} - \overline{y}_k
        \bigr)^\top
        .
    \end{align*}
    The Kalman gain is then
    \begin{align*}
        K_k
        &=
        P^{xy}_k
        \bigl(
            P^{yy}_k
        \bigr)^{-1}.
    \end{align*}
    Finally, the updated ensemble is obtained by shifting each forecast member according to
    \begin{align*}
        \widetilde{X}_k^{(i)}
        &=
        \widehat{X}_k^{(i)}
        +
        K_k
        \bigl(
            O_k - \widehat{Y}_k^{(i)}
        \bigr),
        \qquad i=1,\dots,M.
    \end{align*}
\end{enumerate}
In expectation, this ensemble has the correct posterior mean and covariance under the linear-Gaussian model assumptions \cite{burgers1998analysis}. Iterating this procedure yields an approximate recursive filter. 

We remark that the stochastic EnKF, outlined here, converges to the exact Kalman filter as $M \to \infty$ in the linear Gaussian setting \cite{evensen2003ensemble}. In nonlinear and non-Gaussian problems, it remains a useful approximation but may underestimate uncertainty or produce biased posteriors \cite{katzfuss2016understanding}. In this work, we use the stochastic EnKF as an approximative Gaussian baseline to benchmark the proposed deep density based filters.

\section{Additional example: Linear spring-mass}
\label{SM: LSM}
%\subsection{Linear spring-mass}
%\label{subsec: lsm}
We now consider a linear spring-mass chain with $r = \frac{d}{2}$ masses connected in series by springs and dampers, disturbed by random forcing modelled by the diffusion term in \eqref{eq: state}. We define the system by the constant diffusion coefficient, $\sigma(x) = I_{d\times d}$, and the drift coefficient $\mu(x) = Ax$, where
\begin{align*}
    A
    &=
    \begin{pmatrix}
        0 && I_{r\times r} \\
        A_{21} && A_{22}
    \end{pmatrix}.
\end{align*}
The stiffness and damping matrices $A_{21},A_{22}\in\R^{r\times r}$ are determined by masses $m_1,\dots,m_r$, spring constants $k_1,\dots,k_{r+1}$, and damping coefficients $c_1,\dots,c_{r+1}$ via the standard tridiagonal/diagonal structure:
\begin{align*}
    (A_{21})_{ii} &= -\frac{k_i+k_{i+1}}{m_i}, &
    (A_{21})_{i,i+1} &= \frac{k_{i+1}}{m_i}, &
    (A_{21})_{i+1,i} &= \frac{k_{i+1}}{m_{i+1}}, \\
    (A_{22})_{ii} &= -\frac{c_i+c_{i+1}}{m_i}, &
    (A_{22})_{ij} &= 0 \quad (i\neq j).
\end{align*}
We use heterogeneous parameters sampled independently for each component, distributed as  $m_i \sim \mathrm{Unif}[0.8,1.2]$, $k_i \sim \mathrm{Unif}[0.8,1.2]$, and $c_i \sim \mathrm{Unif}[0.15,0.25]$. This sampling is done once, before running the experiments, and should afterwards be seen as deterministic constants. We let $d'=r$ and define the observation process through the measurement function $h(x)=H x$ with $H=[I_{r\times r}\;\;0_{r\times r}]$, that is, relative positions are observed but not velocities. The initial distribution is $\pi_0 = \N(0,I_{d\times d})$. We fix $T=1$, $K=10$, and $R=I_{r \times r}$. This defines a linear Gaussian state-space model, where the solution $S$ is an Ornstein--Uhlenbeck process, and the solution to the filtering equations is provided by the Kalman filter.

We report results for $d\in\{10,100\}$ in Figure~\ref{fig: LSM}. In this setting we employ the EnKF with $10^5$ ensembles, which in some metrics matches the performance of the LogBSDEF and in others outperforms it. For the $10$-dimensional problem we see decent performance for the LogDSF but evidently the LogBSDEF outperforms it across all metrics, and in the $100$-dimensional case the optimization for the LogDSF does not converge properly and is thus omitted. %In summary, LogBSDEF performs robustly across a range of linear high-dimensional examples. Next, we benchmark the filters in strongly nonlinear settings.

\begin{figure}
    \centering
    \begin{minipage}[t]{0.9\linewidth}  
    \raggedright
    
    \input{include/figures/LSM/groupplot}
    \end{minipage}%
    \begin{minipage}[t]{0.1\linewidth}  
    \vspace{0pt}
    \raggedright
    \begin{tikzpicture}[overlay]
    \node at (0,0) {$ $};
    \node at (-5,-0.5) {\ref*{legend: LSM10d MAE group}};
    \end{tikzpicture}
    \end{minipage}
    \vspace{1.5em}
    \caption{On the left and right panels the results for the $10$-dimensional and $100$-dimensional linear spring-mass models are depicted respectively. From top to bottom the rMAE, FME, and KLD metrics are illustrated.}
    \vspace{-1em}
    \label{fig: LSM}
\end{figure}

\section{Normalization}
\label{SM: Normalization}
We use three different normalization techniques. The first one, quadrature, is only suitable for the one-dimensional examples but is also quite effective in this setting. By defining a connected domain $D = [l,r] \subset \R$ so that
\begin{align*}
    \mathbb{P}
    \big(
        S_t \in D,
        ~t \in [0,T]
    \big) 
    >
    1 - \epsilon,
\end{align*}
for a sufficiently small $\epsilon$, we employ simple Riemann quadrature, based on evenly spaced quadrature points $x^{(i)} \in D,~i=1,\dots,I$. This leads to the following quadrature approximation of the normalization constant:
\begin{align*}
    \widehat{Z}_k^{(m)}
    &=
    \frac{|r-l|}{I}
    \sum_{i=1}^{I}
    \widehat{p}_k
    \bigl(
        x^{(i)} \mid O_{1:k}^{(m)}
    \bigr).
\end{align*}
For higher-dimensional problems, quadrature quickly becomes infeasible due to the curse of dimensionality. In this case, we use importance sampling at each observation time. Given a proposal distribution $q$ and $I \in \mathbb{N}$ proposal samples $\{x^{(i)}\}_{i=1}^{I}$ from $q$, we compute the Monte Carlo estimator
\begin{align*}
    \widehat{Z}_k^{(m)}
    &=
    \frac{1}{I}
    \sum_{i=1}^{I} 
    \frac{
    \widehat{p}_k \big(x^{(i)} \mid O_{1:k}^{(m)} \big)}
    {q\big( x^{(i)} \mid O_{1:k}^{(m)} \big)}
    .
\end{align*}
We consider two proposal choices:
\begin{itemize}
    \item \textbf{EKF-based proposal:} 
    We approximate the filtering distribution with an extended Kalman filter, using its mean and covariance for a Gaussian proposal $q$. This results in observation dependent proposals and requires drawing a fresh set of samples for each observation chain, typically yielding low-variance estimates. 
    
    \item \textbf{Wide-tailed Gaussian proposal:} 
    As a cheaper yet more robust alternative, we use a Gaussian proposal with inflated variance and mean equal to the unconditional mean of $S_{t_k}$ (approximated once from simulated trajectories). Since this proposal is independent of the observation sequence, the same set of proposal samples can be reused for all observation chains, 
    significantly reducing computational cost. However, it provides less relevant samples than the EKF-based proposal.
\end{itemize}

\section{Model architectures}
\label{SM: Architectures}
\begin{figure}[h]
    \centering
    \input{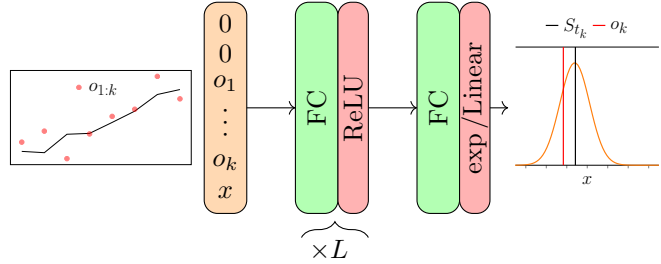}
    \caption{The standard FCN architecture used in the implementation. In the figure, $x$ denotes the state value in which the density is evaluated and $o_{1:k}$ denotes the available observations. The input is padded with zeros so that it has a constant dimension with respect to $k$.}
    \vspace{0em}
    \label{fig: models_fcn}
\end{figure}

To parameterize $(\phi_{k, n})_{k=0, n=1}^{K-1,N}$ in \eqref{eq: deep splitting minimization} for the DSF and $\big(\phi_k,(\overline{v}_{k, n})_{n=0}^{N-1}\big)_{k=0}^{K-1}$ in \eqref{eq: Deep BSDE minimization} for the BSDEF we utilize fully connected neural network architectures. In both cases, $\phi_k \colon \R^{d'\times k} \times \R^d \to \R$ represents a density. Thus we utilize the same architecture, visualized in Figure~\ref{fig: models_fcn}, for both methods. We remark that for all time steps, regardless of the number of available observations, the input dimension is held constant and equal to $d + d' (K-1)$. This allows us to initialize each network with the same weights as the one trained in the previous step, improving convergence and stability across all experiments. To this end, the observations are concatenated with the state value $x$ and remaining unavailable observations are set to 0. The resulting network is rather standard, with $L$ fully connected hidden layers of constant width and the Rectified Linear Unit (ReLU) activation. The final layer has output dimension 1 and for the logarithmic variants we use linear final activation, while the standard DSF and BSDEF have exponential activation to ensure positive output. For BSDEF and its logarithmic version, the networks $(\overline{v}_{k, n})_{n=0}^{N-1}$ have a similar architecture, where the main difference is the output layer, which has dimension $d$ and always linear activation. For these networks representing gradients, we also typically set a lower layer width.

\begin{figure}[h]
    \centering
    \input{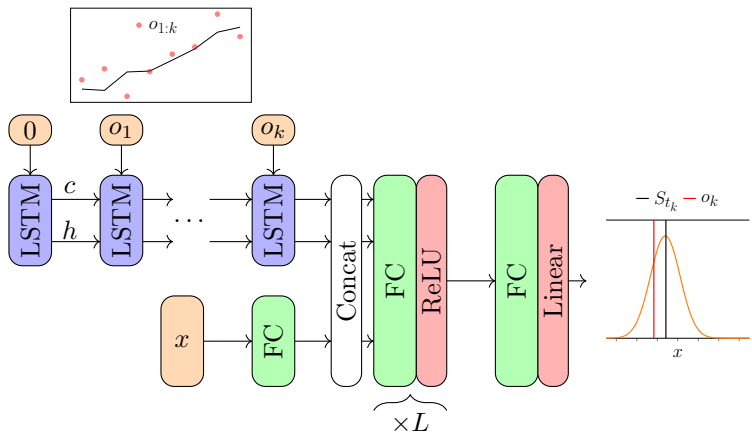}
    \caption{The LSTM-based architecture tested in the long-horizon Ornstein--Uhlenbeck example in 10 dimensions. A token zero input $0 \in \R^{d'}$ is used in the first LSTM cell, while subsequent ones take as input the observation chain $o_{1:k}$.}
    \label{fig: models_lstm}
\end{figure}

In the Ornstein--Uhlenbeck example with 100 observations, we also test a model more suited for capturing long-horizon temporal dependencies, visualized in Figure~\ref{fig: models_lstm}. The model first encodes the observation sequence (started with a token of zeros) with an LSTM network \cite{Hochreiter} and then concatenates the final cell state $c$ and hidden state $h$ with an embedding of the point $x$. This intermediate representation is then decoded using a standard FCN. Since we only apply the LSTM-based architecture for the LogBSDEF, the output layer has linear activation and either $1$ or $d$ neurons depending on whether the model represents $\phi$ or $v$. Note that the LSTM encoder is independent of $x$, which means it only has to be run once when evaluating multiple state values for the same observation sequence. Moreover, to prevent an explosion of the number of parameters, all LSTM encoders for $(v_{k,n})_{n=0}^{N-1}$, with a fixed $k$, share weights.

\section{Training}
\label{SM: Training}

\subsection{Training the deep density methods}
For DSF and LogDSF, we train in epochs over a fixed dataset: we pre-generate $20\,000$ mini-batches and iterate over them for a chosen number of epochs. 
Because the first step after an update typically requires a larger adjustment, we use a longer schedule for the update network (commonly $100$ epochs) and a shorter one for the prediction network (commonly $10$ epochs). This is effective since each newly introduced network is initialized from the previous one, see Supplementary~Materials~\ref{SM: Architectures}.

For BSDEF and LogBSDEF, we instead generate each mini-batch on the fly. This reduces the statistical error compared to reusing the same data but is more computationally demanding. In the limit of sufficiently many samples, both approaches are comparable; however, on-the-fly sampling is particularly costly for DSF/LogDSF because they optimize $N$ objective terms per observation, compared to a single objective for BSDEF/LogBSDEF under identical time discretizations.

A practical consideration is that DSF/LogDSF often use fewer time steps $N$ than BSDE-based models, which we also found to produce more stable results. Throughout this work we choose hyperparameters per (model, example) to achieve the best accuracy. The reported settings below reflect those choices rather than an exhaustive sweep.

This is especially important to note when considering the omitted failed runs of the DSF and LogDSF for some of the high-dimensional or nonlinear examples. Exhaustive attempts were performed to find hyperparameters yielding convergence of the method. However, in some cases we did not find such parameters and this should be further investigated.

We employ early stopping by monitoring the loss (averaged per 200 iterations) and stop training if it does not decrease with a patience of $P$ evaluations. More precisely, we start a counter $p=0$ that increments by one every time the loss increases, and sets it to $p=0$ when the loss decreases again. Coupled to this is the learning rate where we either use a constant rate or a cosine-annealing rate. Our application of the cosine-annealing consist of decreasing the learning rate from $\eta_{\max}$ to $\eta_{\min}$ over a cycle of length $C$ following a cosine-curve:
\begin{align*}
    \eta_c
    =
    \eta_{\min}
    +
    \frac{1}{2}
    \big(
    \eta_{\max}
    -
    \eta_{\min}
    \big)
    \Big(
    1 
    +
    \cos
    \big(\pi \frac{c}{C}\big)
    \Big),
\end{align*}
where $c$ increases every iteration when $p>\frac{P}{2}$. The hyperparameters are set to $P = 50$ and $C = 80$. 

\subsection{Training parameters}
Here we report the parameters affecting the training of the models, used for each model and experiment employed in Section~\ref{Section: Experiments} and Supplementary~Materials~\ref{SM: LSM}. We recall that the
BSDEF-type models (BSDEF, LogBSDEF and its LSTM variant) have one $\phi$-network and $N$ $v$-networks, while the DSF-type models (DSF, LogDSF) only have one $\phi$-network. In Table~\ref{table: training} we report the widths of the networks by Width$_\phi$ and Width$_v$, respectively. In all the examples except for the LSTM networks we use three hidden layers. We abbreviate the normalization methods of Supplementary~Materials~\ref{SM: Normalization} with I-EKF and I-G for the EKF-based and wide-tailed Gaussian proposals respectively.

Regarding the LSTM-based model in the long-horizon OU example of Section~\ref{section: longrange highdim OU}, it has $4$ hidden layers in the FCN decoder with constant width 512 for the $\phi$-network and 128 for the $v$-networks. The embedding layer for the state value $x$ has the same width as the hidden layers in the decoder and the LSTM encoder has hidden dimension $256$ and $128$ in the $\phi$- and $v$-cases respectively. The LSTM encoder also has 2 extra cells without any input between each observation cell. The other hyperparameters ($N$, $B$, LR, Normalization method) are the same as the FCN for the same example, see Table \ref{table: training}.

\begin{table}[h!]
\centering
\caption{Training setup across models and examples. Models are abbreviated B - BSDEF, LB - LogBSDEF, D - DSF, LD - LogDSF. Columns: time steps $N$, widths ($W_\phi$, $W_v$) for the $\phi$- and $v$-networks (3 hidden layers fixed), batch size $B$, LR schedule (start $\to$ end), and normalization (method and number of normalized trajectories per batch).}
\label{table: training}
\setlength{\tabcolsep}{4.5pt}
\begin{tabular}{l l c c c c c c}
\toprule
\textbf{Model} & \textbf{Example} & $\boldsymbol{N}$ & \textbf{$W_\phi$} & \textbf{$W_v$} &  $\boldsymbol{B}$ & \textbf{LR schedule} & \textbf{Norm} \\
\midrule
\multicolumn{7}{l}{\emph{BSDEF models}} \\
\midrule
B         & OU (1d)               & 64 & 128 & 32 & 512   & Const $10^{-4}$ & (Quad, 64) \\
B         & Bistable         & 64 & 128 & 32 & 512   & Const $10^{-4}$ & (Quad, 64) \\
\midrule
LB      & OU (1d)               & 64 & 128 & 32 & 512   & Const $10^{-4}$ & (Quad, 64) \\
LB      & Bistable         & 64 & 128 & 32 & 512   & Const $10^{-4}$ & (Quad, 64) \\
LB      & Schlögl               & 128 & 256 & 64 & 4096   & Cos $10^{-4}\!\to\!10^{-6}$ & (Quad, 128) \\
\midrule
LB      & OU (10d)              & 16 & 512 & 256 &  1024  & Cos $10^{-4}\!\to\!10^{-5}$ & (I-EKF, 64) \\
LB      & OU (100d)             & 32 & 1024 & 512 &  2048  & Cos $10^{-4}\!\to\!10^{-5}$ & (I-EKF, 64) \\
LB      & LSM (10d)             & 32 & 512 & 256 &  2048  & Cos $10^{-4}\!\to\!10^{-5}$ & (I-EKF, 64) \\
LB      & LSM (100d)            & 32 & 2048 & 512 & 4096   & Cos $10^{-4}\!\to\!10^{-5}$ & (I-EKF, 512) \\
\midrule
LB      & L96 (4d)       & 32 & 256 & 64 & 2048   & Const $10^{-4}$ & (I-G, 64) \\
LB      & L96 (10d)      & 32 & 512 & 256 & 2048   & Const $10^{-4}$ & (I-G, 64) \\
LB      & L96 (20d)      & 32 & 512 & 256 & 2048   & Const $10^{-4}$ & (I-G, 64) \\
LB      & L96 (40d)      & 32 & 512 & 256 & 2048   & Const $10^{-4}$ & (I-G, 64) \\
LB      & L96 (100d)     & 32 & 1024 & 512 & 2048   & Const $10^{-4}$ & (I-G, 64) \\
\midrule
\multicolumn{7}{l}{\emph{DSF models}} \\
\midrule
D           & OU (1d)               & 16               & 128 & N/A             & 2048  & Const $10^{-4}$ & (Quad, 64) \\
D           & Bistable         & 16               & 128 & N/A             & 2048  & Const $10^{-4}$ & (Quad, 64) \\
\midrule
LD           & OU (1d)               & 16               & 128 & N/A             & 2048  & Const $10^{-4}$ & (Quad, 64) \\
LD           & Bistable         & 16               & 128 & N/A             & 2048  & Const $10^{-4}$ & (Quad, 64) \\
LD        & OU (10d)              & 16                & 512  & N/A             & 1024   & Const $10^{-4}$ & (I-EKF, 64) \\
LD        & LSM (10d)             & 32                & 512  & N/A             & 2048    & Const $10^{-4}$ & (I-EKF, 64) \\
\bottomrule
\end{tabular}
\end{table}

\section{Evaluation}
\label{SM: Evaluation}
In this section we briefly detail the parameters used for evaluating the metrics of Section~\ref{section: experiment setup}. To evaluate the probability density of the particle filters and ensemble Kalman filters, both for the benchmark methods and the reference solution, we employ Gaussian kernel density estimators, implemented using \texttt{scipy.stats.gaussian\_kde}~\cite{2020SciPy-NMeth-short}. This function applies Scott’s rule for bandwidth selection by default \cite{Scott1992, silverman1986density}, and we perform no additional tuning of the bandwidth parameter. In Table~\ref{table: evaluation setup} we report the evaluation parameters used for all the experiments. 

\begin{table}[h!]
    \centering
    \caption{Evaluation setup across examples. Columns: number of Monte Carlo evaluation samples $M$, normalization method, reference solution, reference discretization steps $N_{\text{ref}}$. For the Lorenz-96 high-dimensional example we do not have access to a reference solution and only evaluate the MAE and NLL metrics.}
    \label{table: evaluation setup}
    \begin{tabular}{lcccc}
        \toprule
        \textbf{Example} 
        & \textbf{$\boldsymbol{M}$} 
        & \textbf{Normalization method} 
        & \textbf{Reference} 
        & $\boldsymbol{N_{\text{ref}}}$ \\
        \midrule
        OU (1d)      
        & $10^{4}$ 
        & Quad
        & KF 
        & 512 \\
        
        Bistable
        & $10^{4}$ 
        & Quad
        & PF $10^{6}$
        & 128  \\
        
        OU (10d)
        & $10^{3}$ 
        & I-EKF
        & KF 
        & 128 \\
        
        OU (100d)
        & $10^{3}$ 
        & I-EKF 
        & KF 
        & 128 \\
        
        LSM (10d)
        & $10^{3}$ 
        & I-EKF
        & KF 
        & 128 \\
        
        LSM (100d)
        & $10^{3}$ 
        & I-EKF
        & KF 
        & 128 \\
        
        Schlögl
        & $10^{4}$ 
        & Quad
        & PF $10^{6}$
        & 128 \\
        
        Lorenz-96 (4d)
        & $10^{3}$ 
        & I-G
        & PF $10^{6}$
        & 128 \\
        
        \midrule
        \multicolumn{5}{l}{\emph{Only MAE and NLL}} \\
        \midrule
        Lorenz-96 (4d)
        & $10^{4}$ 
        & I-G
        & N/A
        & N/A \\
        
        Lorenz-96 (10d)
        & $10^{4}$ 
        & I-G
        & N/A
        & N/A \\

        Lorenz-96 (20d)
        & $10^{4}$ 
        & I-G
        & N/A
        & N/A \\

        Lorenz-96 (40d)
        & $10^{4}$ 
        & I-G
        & N/A
        & N/A \\

        Lorenz-96 (100d)
        & $10^{4}$ 
        & I-G
        & N/A
        & N/A \\
        \bottomrule
    \end{tabular}
\end{table}

\section{Additional figures}
\label{SM: densities}
\begin{figure}
    \centering
    \input{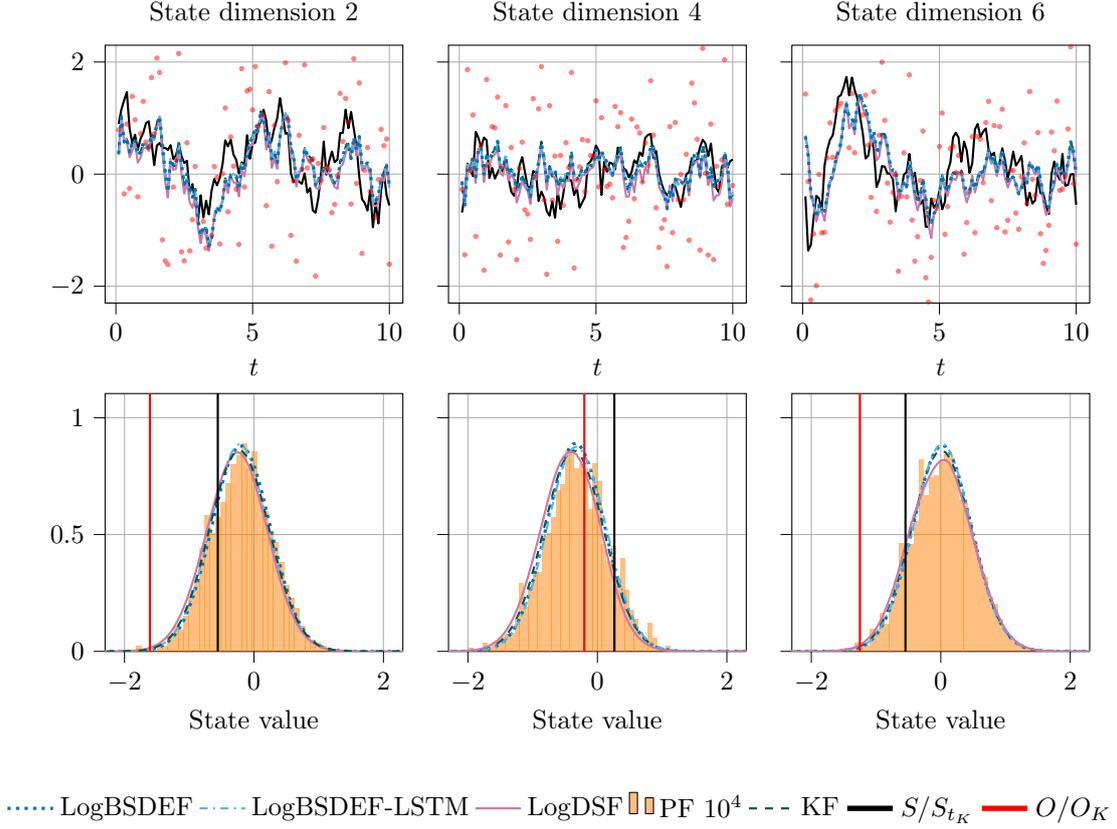}
    \vspace{0em}
    \caption{One trajectory for the long-horizon Ornstein--Uhlenbeck problem in 10 dimensions. The left column shows the sample path with corresponding filter mean estimates. The right column displays the filtering densities at the final time $T=10$.}
    \vspace{-1em}
    \label{fig: OU10d_traj}
\end{figure}

\begin{figure}
    \centering
    \input{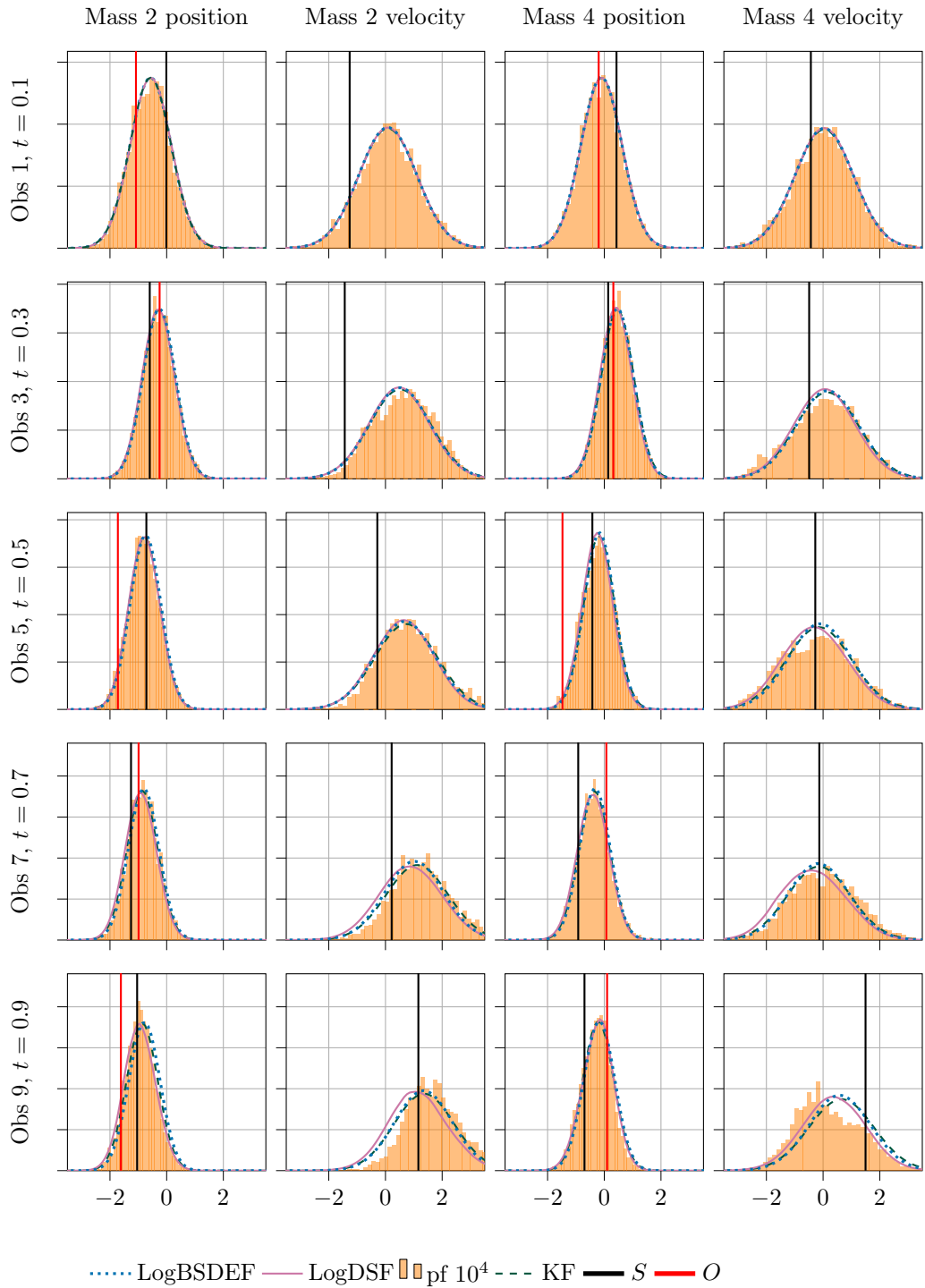}
    \vspace{-1em}
    \caption{Marginal densities over time in the $10$-dimensional linear spring-mass example for selected components. Recall that only positions are observed.}
    \label{fig:pdfs_LSM_10}
\end{figure}

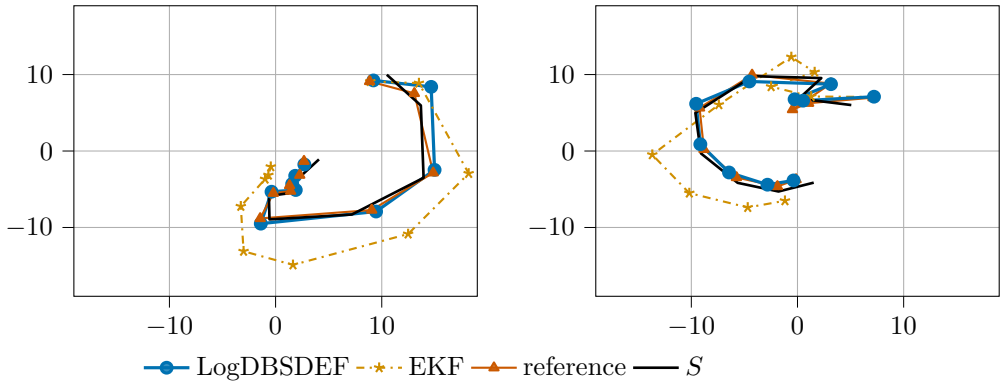
\begin{figure}[h]
\captionsetup{justification=raggedleft}
\noindent
\begin{subfigure}[b]{0.45\linewidth}
\centering
% This file was created with tikzplotlib v0.10.1.
\begin{tikzpicture}

\definecolor{darkgray176}{RGB}{176,176,176}
\definecolor{green}{RGB}{0,128,0}
\definecolor{purple}{RGB}{128,0,128}

\begin{axis}[
width=0.40*6.028in,
height=0.40*4.754in,
legend cell align={left},
legend style={
  fill opacity=0.8,
  draw opacity=1,
  text opacity=1,
  at={(0.03,0.97)},
  anchor=north west,
  draw=none,
  nodes={scale=1.0, transform shape},
  legend columns=-1
},
tick align=outside,
tick pos=left,
x grid style={darkgray176},
xmajorgrids,
ymajorgrids,
xmin=-19, xmax=19,
xtick style={color=black},
y grid style={darkgray176},
ymin=-19, ymax=19,
ytick style={color=black},
legend to name={legend: Lorenz 4 paths},
]
\addplot[method LogBSDEF]
table {%
9.21315237994671 9.23208512418138
14.6540954246544 8.39123132292781
14.9748023938755 -2.46255071348022
9.46352523312385 -7.91484717096626
-1.39271093437353 -9.50210631715903
-0.381868570391086 -5.32364330867907
1.8997128785439 -5.09958313795162
1.5643909592128 -4.29956816544126
1.84968617174117 -3.23957923718534
2.70937028411326 -1.76223458930743
};
\addlegendentry{LogDBSDEF}
%\addplot [semithick, blue, mark=*, mark size=1, mark options={solid}]
\addplot [method EKF]
table {%
8.91772784740391 9.08956558486106
13.5018409907904 8.88365353528359
18.185708989957 -2.935375750804
12.4830806603354 -10.8541924133632
1.63789030201378 -14.8931223684487
-3.01608824910497 -13.1118825135883
-3.27493028600388 -7.23383931277065
-0.993006072719553 -3.67276994496215
-0.677737144812872 -3.15956845490838
-0.453232899889304 -2.05710428135596
};
\addlegendentry{EKF}
%\addplot [semithick, green, mark=*, mark size=1, mark options={solid}]
\addplot[method PF 1e6]
table {%
8.85253853854789 9.10395542344274
13.0719525714135 7.54196384395089
14.8431739714419 -2.84484433988404
9.05444046936931 -7.71905910242362
-1.45030113746901 -8.81721572936819
-0.209105185631429 -5.51462592482994
1.41489720805158 -5.23857470592105
1.31746248373702 -4.51039103132697
2.2745271224219 -3.15336282321127
2.67684501731899 -1.33242105640204
};
\addlegendentry{reference}
%\addplot [semithick, red, mark=*, mark size=1, mark options={solid}]
\addplot[line width=1pt, draw=black]
table {%
10.4771560410025 9.97807352608754
13.6825250344828 5.96028614094973
13.9411046344057 -3.50020951103188
7.18145055750569 -8.30706845618564
-0.574945746098028 -8.94677045393922
-0.626469977545929 -5.87320493961675
1.24328363494185 -5.52527040825402
1.44416826377209 -4.80399739445122
2.13287901549163 -3.5863643826944
4.08607465930433 -1.09453492917377
};
\addlegendentry{$S$}
\end{axis}

\end{tikzpicture}
\captionsetup{margin=90pt}
\vspace{0pt}
\end{subfigure}
\begin{subfigure}[b]{0.45\linewidth}
\centering
% This file was created with tikzplotlib v0.10.1.
\begin{tikzpicture}

\definecolor{darkgray176}{RGB}{176,176,176}
\definecolor{green}{RGB}{0,128,0}
\definecolor{purple}{RGB}{128,0,128}

\begin{axis}[
width=0.40*6.028in,
height=0.40*4.754in,
legend cell align={left},
legend style={
  fill opacity=0.8,
  draw opacity=1,
  text opacity=1,
  at={(0.03,0.97)},
  anchor=north west,
  draw=none,
  nodes={scale=1.0, transform shape},
  legend columns=-1
},
tick align=outside,
tick pos=left,
xmajorgrids,
ymajorgrids,
x grid style={darkgray176},
xmin=-19, xmax=19,
xtick style={color=black},
y grid style={darkgray176},
ymin=-19, ymax=19,
ytick style={color=black},
legend to name={legend: not used},
]
%\addplot [semithick, blue, mark=*, mark size=1, mark options={solid}]
\addplot[method EKF]
table {%
7.30018295634956 7.07242257591924
1.28778301615874 7.12896380919914
-2.48362944524058 8.41677710474459
1.61704845957486 10.3218612351324
-0.580239688122597 12.2966153675204
-7.40908632822367 6.05897721382446
-13.676232439423 -0.483469214295791
-10.1884410982658 -5.5015308642945
-4.68646844764276 -7.38422038847951
-1.18065582947254 -6.50439963015119
};
\addlegendentry{EKF}
%\addplot [semithick, green, mark=*, mark size=1, mark options={solid}]
\addplot[method PF 1e6]
table {%
7.10610570725643 6.98093798386397
1.0495806439029 6.25437633066317
-0.470883384896558 5.42745224529671
3.14599778079307 8.82097746256948
-4.27615506382307 9.94637502382057
-9.22495506187112 5.64176983385838
-8.81295097387948 0.223042738632962
-5.74913142829406 -3.46061457624128
-1.88660383737213 -4.63554754121081
-0.133192246375716 -3.73040568141125
};
\addlegendentry{Reference}
%\addplot [semithick, purple, mark=*, mark size=1, mark options={solid}]
\addplot[method LogBSDEF]
table {%
7.20502609039301 7.10334021242894
0.528285056747708 6.58541258555774
-0.270645085181149 6.7709144662643
3.15474980385458 8.74117866893254
-4.53030791069472 9.097110693721
-9.54484508549891 6.17850799376514
-9.15046664764773 0.885295862418788
-6.44724624890211 -2.80766947274349
-2.82598744000192 -4.40497954416698
-0.370799757105526 -3.84799714522787
};
\addlegendentry{LogDBSDEF}
%\addplot [semithick, red, mark=*, mark size=1, mark options={solid}]
\addplot[line width=1pt, draw=black]
table {%
5.06888632832034 6.01066279624782
0.630777618375558 6.77293221996659
0.123638454131753 6.83910813002036
2.24457281213631 9.54450789046062
-4.22752764408406 9.79020463277612
-9.61881361137033 4.98744326388194
-9.09062840241241 -0.313177923719736
-5.63562583763802 -4.16217278503704
-1.78592775029414 -5.29392931774134
1.5289902977663 -4.13678574206661
};
\addlegendentry{$S$}
\end{axis}

\end{tikzpicture}
\captionsetup{margin=80pt}
\vspace{0pt}
\end{subfigure}
\begin{tikzpicture}[overlay]
\node at (0,0) {$ $};
\node at (-8.0,-0.3) {\ref*{legend: Lorenz 4 paths}};
\end{tikzpicture}
\vspace{10pt}
\captionsetup{justification=justified}
\caption{One realization of the state process $S$, starting at $S_0\sim \N\big((8,8,8,8),I\big)$, in the four-dimensional Lorenz-96 example together with the corresponding filtering means from the different methods. The left and right panels illustrate the first two and last two variables respectively.}
\label{fig: Lorenz 4 paths}
\end{figure}

In this section we present qualitative results: example trajectories and density comparisons. Figure~\ref{fig: OU10d_traj} shows a trajectory and the final-time density for the long-horizon $10$-dimensional Ornstein--Uhlenbeck model with some different methods. Figures~\ref{fig:LSM_traj} and \ref{fig:pdfs_LSM_10} analogously show trajectories and densities for the $10$-dimensional linear spring-mass model. Finally, Figure~\ref{fig: Lorenz 4 paths} shows one trajectory of the four-dimensional Lorenz-96 model in the $x_1 x_2$- and $x_3 x_4$-marginal planes.

\begin{figure}
    \centering
    \input{include/figures/LSM10/trajectories}
    \vspace{-0.1cm}
    \caption{One trajectory and filter estimates for selected components in the $10$-dimensional linear spring-mass example.}
    \vspace{-0.5cm}
    \label{fig:LSM_traj}
\end{figure}

\section{The chemical Schl\"ogl model}
\label{SM: schlogl chemistry}
The Schlögl model \cite{schlogl1972chemical, vellela2009stochastic, Vlysidis2018} is a prototypical bistable chemical reaction network involving a single reacting species $S$ coupled to two external reservoirs $A$ and $B$. Its reactions are
\begin{align}
    \label{eq: schlogl reaction one}
        2S + A &\xrightarrow{\theta_1} 
        3S,
        & 3S 
        &\xrightarrow{\theta_2} 
        2S + A, 
        \\
    \label{eq: schlogl reaction two}
        B 
        &\xrightarrow{\theta_3} 
        S,
        & 
        S 
        &\xrightarrow{\theta_4} 
        B.
\end{align}
Here, the notation $X \xrightarrow{\theta_n} Y$ means that configuration $X$ is transformed into configuration $Y$ with rate constant $\theta_n$ according to mass-action kinetics.

Reactions in \eqref{eq: schlogl reaction one} describe an autocatalytic feedback loop: two molecules of $S$ react with a molecule from reservoir $A$ to produce a third $S$, while the reverse reaction consumes three $S$ to regenerate $A$. Reactions in \eqref{eq: schlogl reaction two} model inflow and outflow of $S$ through coupling with reservoir $B$, maintaining a constant supply and removal of the species.

Starting from the chemical master equation and applying the chemical Langevin approximation \cite{gillespie2000chemical}, one obtains the SDE, given by the dynamics detailed in Section~\ref{section: schlogl}, governing the concentration of $S$. Each reaction channel contributes an independent noise source $(B^{(i)})_{i=1}^4$, corresponding to intrinsic stochastic fluctuations of the reaction events. The nonlinear drift encodes autocatalytic growth and saturation, while the multiplicative diffusion terms capture random reaction noise.

For parameter values
\begin{align*}
    \theta 
    =
    (3\times 10^{-7},\, 10^{-4},\, 10^{-3},\, 3.5),
    \quad
    A = 10^{5},
    \quad
    B = 2\times 10^{5},
\end{align*}
the deterministic mean-field limit exhibits two stable steady states, separated by an unstable equilibrium. Random fluctuations can drive rare transitions between these states, a hallmark of bistable systems.  
This bistability makes the Schlögl model a canonical benchmark in stochastic chemical kinetics and nonlinear filtering of reaction networks.

%\putbib
%\end{bibunit}

\printbibliography[title={References}]

@article{andersson2026deep,
  title={The deep multi-{FBSDE} method: a robust deep learning method for coupled {FBSDEs}},
  author={Andersson, Kristoffer and Andersson, Adam and Oosterlee, Cornelis W},
  journal={J. Sci. Comput.},
  volume={106},
  number={3},
  pages={77},
  year={2026},
  publisher={Springer}
}

@article {bao2024score,
    AUTHOR = {Bao, Feng and Zhang, Zezhong and Zhang, Guannan},
     TITLE = {A score-based filter for nonlinear data assimilation},
   JOURNAL = {J. Comput. Phys.},
  FJOURNAL = {Journal of Computational Physics},
    VOLUME = {514},
      YEAR = {2024},
     PAGES = {Paper No. 113207, 16},
      noISSN = {0021-9991,1090-2716},
   MRCLASS = {65C35 (62M20 65C05)},
  MRNUMBER = {4762023},
       noDOI = {10.1016/j.jcp.2024.113207},
       noURL = {https://doi.org/10.1016/j.jcp.2024.113207},
}

@book{barshalom2001estimation,
  author    = {Yaakov Bar-Shalom and X. Rong Li and Thiagalingam Kirubarajan},
  title     = {Estimation with Applications to Tracking and Navigation},
  publisher = {John Wiley \& Sons},
  year      = {2001}
}

@article {beck2019machine,
    AUTHOR = {Beck, Christian and E, Weinan and Jentzen, Arnulf},
     TITLE = {Machine learning approximation algorithms for high-dimensional
              fully nonlinear partial differential equations and
              second-order backward stochastic differential equations},
   JOURNAL = {J. Nonlinear Sci.},
  FJOURNAL = {Journal of Nonlinear Science},
    VOLUME = {29},
      YEAR = {2019},
    NUMBER = {4},
     PAGES = {1563--1619},
      noISSN = {0938-8974,1432-1467},
   MRCLASS = {65Z05 (35K55 35R60 60H10)},
  MRNUMBER = {3993178},
       noDOI = {10.1007/s00332-018-9525-3},
       noURL = {https://doi.org/10.1007/s00332-018-9525-3},
}

@article {brajard2020combining,
    AUTHOR = {Brajard, Julien and Carrassi, Alberto and Bocquet, Marc and
              Bertino, Laurent},
     TITLE = {Combining data assimilation and machine learning to emulate a
              dynamical model from sparse and noisy observations: a case
              study with the {L}orenz 96 model},
   JOURNAL = {J. Comput. Sci.},
  FJOURNAL = {Journal of Computational Science},
    VOLUME = {44},
      YEAR = {2020},
     PAGES = {101171, 11},
      noISSN = {1877-7503,1877-7511},
   MRCLASS = {86A22 (37M10 62M20)},
  MRNUMBER = {4117875},
       noDOI = {10.1016/j.jocs.2020.101171},
       noURL = {https://doi.org/10.1016/j.jocs.2020.101171},
}

@article{bagmark_1,
  title={An energy-based deep splitting method for the nonlinear filtering problem},
  author={B{\aa}gmark, Kasper and Andersson, Adam and Larsson, Stig},
  journal={Partial Differ. Equ. Appl.},
  volume={4},
  nonumber={2},
  nopages={14},
  year={2023},
  publisher={Springer}
}

@article{baagmark2024convergent,
  title={A convergent scheme for the {B}ayesian filtering problem based on the {F}okker--{P}lanck equation and deep splitting},
  author={B{\aa}gmark, Kasper and Andersson, Adam and Larsson, Stig and Rydin, Filip},
  journal={arXiv:2409.14585},
  year={2024}
}

@article{baagmark2025nonlinear,
  title={Nonlinear filtering based on density approximation and deep {BSDE} prediction},
  author={B{\aa}gmark, Kasper and Andersson, Adam and Larsson, Stig},
  journal={arXiv:2508.10630},
  year={2025}
}

@article{Arnulf,
  title={Deep learning based numerical approximation algorithms for stochastic partial differential equations and high-dimensional nonlinear filtering problems},
  author={Beck, Christian and Becker, Sebastian and Cheridito, Patrick and Jentzen, Arnulf and Neufeld, Ariel},
  journal={arXiv:2012.01194},
  year={2020}
}

@article{Arnulf_PDE,
  title={Deep splitting method for parabolic {PDE}s},
  author={Beck, Christian and Becker, Sebastian and Cheridito, Patrick and Jentzen, Arnulf and Neufeld, Ariel},
  journal={SIAM J. Sci. Comput.},
  volume={43},
  nonumber={5},
  pages={A3135--A3154},
  year={2021},
  publisher={SIAM}
}

@incollection {bickel2008sharp,
    AUTHOR = {Bickel, Peter and Li, Bo and Bengtsson, Thomas},
     TITLE = {Sharp failure rates for the bootstrap particle filter in high
              dimensions},
 BOOKTITLE = {Pushing the limits of contemporary statistics: contributions
              in honor of {J}ayanta {K}. {G}hosh},
    SERIES = {Inst. Math. Stat. (IMS) Collect.},
    VOLUME = {3},
     PAGES = {318--329},
 PUBLISHER = {Inst. Math. Statist., Beachwood, OH},
      YEAR = {2008},
      noISBN = {978-0-940600-75-1},
   MRCLASS = {93E11 (60G50 62L12)},
  MRNUMBER = {2459233},
       noDOI = {10.1214/074921708000000228},
       noURL = {https://doi.org/10.1214/074921708000000228},
}

@book{blackman1999design,
  title={Design and Analysis of Modern Tracking Systems},
  author={Blackman, Samuel S and Popoli, Robert},
  year={1999},
  publisher={Artech House Publishers}
}

@article {burgers1998analysis,
      author = "Gerrit Burgers and Peter Jan van Leeuwen and Geir Evensen",
      title = "Analysis scheme in the ensemble {K}alman filter",
      journal = "Mon. Wea. Rev.",
      year = "1998",
      publisher = "American Meteorological Society",
      address = "Boston MA, USA",
      volume = "126",
      number = "6",
      nodoi = "10.1175/1520-0493(1998)126<1719:ASITEK>2.0.CO;2",
      pages=      "1719 - 1724",
      nourl = "https://journals.ametsoc.org/view/journals/mwre/126/6/1520-0493_1998_126_1719_asitek_2.0.co_2.xml"
}

@article{chan2019machine,
  title={Machine learning for semi linear {PDE}s},
  author={Chan-Wai-Nam, Quentin and Mikael, Joseph and Warin, Xavier},
  journal={J. Sci. Comput.},
  volume={79},
  number={3},
  pages={1667--1712},
  year={2019},
  publisher={Springer}
}

@article {chopin2004central,
    AUTHOR = {Chopin, Nicolas},
     TITLE = {Central limit theorem for sequential {M}onte {C}arlo methods
              and its application to {B}ayesian inference},
   JOURNAL = {Ann. Statist.},
  FJOURNAL = {The Annals of Statistics},
    VOLUME = {32},
      YEAR = {2004},
    NUMBER = {6},
     PAGES = {2385--2411},
      noISSN = {0090-5364,2168-8966},
   MRCLASS = {60F05 (62F15 65C05)},
  MRNUMBER = {2153989},
MRREVIEWER = {C.\ S.\ Withers},
       noDOI = {10.1214/009053604000000698},
       noURL = {https://doi.org/10.1214/009053604000000698},
}

@book{CoverThomas2006,
  author    = {Thomas M. Cover and Joy A. Thomas},
  title     = {Elements of Information Theory},
  edition   = {2nd},
  publisher = {Wiley-Interscience},
  year      = {2006},
  noISBN      = {978-0-471-24195-9},
  nodoi       = {10.1002/047174882X},
  nourl       = {https://doi.org/10.1002/047174882X}
}

@article{cui2005comparison,
  title={A comparison of nonlinear filtering approaches with an application to ground target tracking},
  author={Cui, Ningzhou and Hong, Lang and Layne, Jeffery R},
  journal={Signal Processing},
  volume={85},
  nonumber={8},
  pages={1469--1492},
  year={2005},
  publisher={Elsevier}
}

@article {crisan2002survey,
    AUTHOR = {Crisan, Dan and Doucet, Arnaud},
     TITLE = {A survey of convergence results on particle filtering methods
              for practitioners},
   JOURNAL = {IEEE Trans. Signal Process.},
  FJOURNAL = {IEEE Transactions on Signal Processing},
    VOLUME = {50},
      YEAR = {2002},
    NUMBER = {3},
     PAGES = {736--746},
      noISSN = {1053-587X,1941-0476},
   MRCLASS = {60F15 (62F15)},
  MRNUMBER = {1895071},
       noDOI = {10.1109/78.984773},
       noURL = {https://doi.org/10.1109/78.984773},
}

@article{Duc_Kuroda,
  title={Ensemble  {Kalman} Filter data assimilation and storm surge experiments of tropical cyclone Nargis},
  author={Duc, Le and Kuroda, Tohru and Saito, Kazuo and Fujita, Tadashi},
  journal={Tellus A},
  volume={67},
  nonumber={1},
  pages={25941},
  year={2015},
  publisher={Taylor \& Francis}
}

@article{E_2017,
   title={Deep Learning-Based Numerical Methods for High-Dimensional Parabolic Partial Differential Equations and Backward Stochastic Differential Equations},
   volume={5},
   noISSN={2194-671X},
   nourl={http://dx.doi.org/10.1007/s40304-017-0117-6},
   noDOI={10.1007/s40304-017-0117-6},
   nonumber={4},
   journal={Commun. Math. Stat},
   publisher={Springer Science and Business Media LLC},
   author={E, Weinan and Han, Jiequn and Jentzen, Arnulf},
   year={2017},
   month=nov, pages={349–380} 
}

@article{ehrendorfer2007review,
  title={A review of issues in ensemble-based {K}alman filtering},
  author={Ehrendorfer, Martin},
  journal={Meteorol. Z.},
  volume={16},
  year={2007}
}

@article{evensen2003ensemble,
  title={The ensemble {K}alman filter: {T}heoretical formulation and practical implementation},
  author={Evensen, Geir},
  journal={Ocean Dyn.},
  volume={53},
  number={4},
  pages={343--367},
  year={2003},
  publisher={Springer}
}

@article{evensen1994sequential,
author = {Evensen, Geir},
title = {Sequential data assimilation with a nonlinear quasi-geostrophic model using {M}onte {C}arlo methods to forecast error statistics},
journal = {J. Geophys. Res.},
volume = {99},
number = {C5},
pages = {10143-10162},
nodoi = {https://doi.org/10.1029/94JC00572},
nourl = {https://agupubs.onlinelibrary.wiley.com/doi/abs/10.1029/94JC00572},
noeprint = {https://agupubs.onlinelibrary.wiley.com/doi/pdf/10.1029/94JC00572},
year = {1994}
}

@article {frey2022convergence,
    AUTHOR = {Frey, R\"udiger and K\"ock, Verena},
     TITLE = {Convergence analysis of the deep splitting scheme: the case of
              partial integro-differential equations and the associated
              forward backward {SDE}s with jumps},
   JOURNAL = {SIAM J. Sci. Comput.},
  FJOURNAL = {SIAM Journal on Scientific Computing},
    VOLUME = {47},
      YEAR = {2025},
    NUMBER = {1},
     PAGES = {A527--A552},
      noISSN = {1064-8275,1095-7197},
   MRCLASS = {65C30 (60H35 65M12)},
  MRNUMBER = {4865059},
       noDOI = {10.1137/23M1595710},
       noURL = {https://doi.org/10.1137/23M1595710},
}

@article{galanis2006applications,
  title={{Applications of Kalman filters based on non-linear functions to numerical weather predictions}},
  author={Galanis, George and Louka, Petroula and Katsafados, P and Pytharoulis, I and Kallos, G},
  journal={Ann. Geophys},
  volume={24},
  pages={1--10},
  year={2006},
  publisher={Citeseer}
}

@article {Germain,
    AUTHOR = {Germain, Maximilien and Pham, Huy\^en and Warin, Xavier},
     TITLE = {Approximation error analysis of some deep backward schemes for
              nonlinear {PDE}s},
   JOURNAL = {SIAM J. Sci. Comput.},
  FJOURNAL = {SIAM Journal on Scientific Computing},
    VOLUME = {44},
      YEAR = {2022},
    NUMBER = {1},
     PAGES = {A28--A56},
      noISSN = {1064-8275,1095-7197},
   MRCLASS = {65M99 (60H35 65C20 65M12)},
  MRNUMBER = {4358471},
MRREVIEWER = {Xiaobing\ Henry\ Feng},
       noDOI = {10.1137/20M1355355},
       noURL = {https://doi.org/10.1137/20M1355355},
}

@article{gillespie2000chemical,
  title={The chemical {L}angevin equation},
  author={Gillespie, Daniel T},
  journal={J. Chem. Phys.},
  volume={113},
  number={1},
  pages={297--306},
  year={2000},
  publisher={American Institute of Physics}
}

@book {goodman1997mathematics,
    AUTHOR = {Goodman, I. R. and Mahler, Ronald P. S. and Nguyen, Hung T.},
     TITLE = {Mathematics of Data Fusion},
    SERIES = {Theory and Decision Library. Series B: Mathematical and Statistical Methods},
    VOLUME = {37},
 PUBLISHER = {Kluwer Academic Publishers Group, Dordrecht},
      YEAR = {1997},
     PAGES = {xii+507},
      noISBN = {0-7923-4674-2},
   MRCLASS = {60D05 (03B52 28E10 62-07 68T35)},
  MRNUMBER = {1635258},
MRREVIEWER = {Giulianella\ Coletti},
       noDOI = {10.1007/978-94-015-8929-1},
       noURL = {https://doi.org/10.1007/978-94-015-8929-1},
}

@article{gordon1993novel,
  title={Novel approach to nonlinear/non-{G}aussian {B}ayesian state estimation},
  author={Gordon, Neil J. and Salmond, David J. and Smith, Adrian F. M.},
  journal={IEEE Proceedings F (Radar and Signal Processing)},
  volume={140},
  number={2},
  pages={107--113},
  year={1993},
  publisher={IET},
  nodoi={10.1049/ip-f-2.1993.0015}
}

@article{HAN20081434,
    title = {An evaluation of the nonlinear/non-{G}aussian filters for the sequential data assimilation},
    journal = {Remote Sens. Environ.},
    volume = {112},
    number = {4},
    pages = {1434-1449},
    year = {2008},
    nonote = {Remote Sensing Data Assimilation Special Issue},
    noISSN = {0034-4257},
    nodoi = {https://doi.org/10.1016/j.rse.2007.07.008},
    nourl = {https://www.sciencedirect.com/science/article/pii/S0034425707003409},
    author = {Xujun Han and Xin Li}
}

@article{han2025brief,
  title={A brief review of the Deep {BSDE} method for solving high-dimensional partial differential equations},
  author={Han, Jiequn and Jentzen, Arnulf and E, Weinan},
  journal={arXiv:2505.17032},
  year={2025}
}

@article {han_convergence,
    AUTHOR = {Han, Jiequn and Long, Jihao},
     TITLE = {Convergence of the deep {BSDE} method for coupled {FBSDE}s},
   JOURNAL = {Probab. Uncertain. Quant. Risk},
  FJOURNAL = {Probability, Uncertainty and Quantitative Risk},
    VOLUME = {5},
      YEAR = {2020},
     PAGES = {Paper No. 5, 33},
      noISSN = {2095-9672,2367-0126},
   MRCLASS = {60H35 (60H10 65C30)},
  MRNUMBER = {4122227},
       noDOI = {10.1186/s41546-020-00047-w},
       noURL = {https://doi.org/10.1186/s41546-020-00047-w},
}

@article{Hochreiter,
    author = {Hochreiter, Sepp and Schmidhuber, Jürgen},
    title = {Long Short-Term Memory},
    journal = {Neural Computation},
    volume = {9},
    number = {8},
    pages = {1735-–1780},
    year = {1997}
}

@incollection{johannes2009mcmc,
  author    = {Michael S. Johannes and Nicholas G. Polson},
  title     = {{MCMC} Methods for Continuous-Time Financial Econometrics},
  booktitle = {Handbook of Financial Econometrics},
  publisher = {Elsevier},
  year      = {2009},
  editor_notused    = {Yacine Aït-Sahalia and Lars P. Hansen},
  pages     = {1--72}
}

@article{kamino2023optimal,
  title={Optimal inference of molecular interaction dynamics in {FRET} microscopy},
  author={Kamino, Keita and Kadakia, Nirag and Avgidis, Fotios and Liu, Zhe-Xuan and Aoki, Kazuhiro and Shimizu, Thomas S and Emonet, Thierry},
  journal={Proc. Natl. Acad. Sci. U.S.A.},
  volume={120},
  number={15},
  eid={e2211807120},
  year={2023},
  publisher={National Academy of Sciences}
}

@article{karimi2010extensive,
  author  = {Karimi, A. and Paul, M. R.},
  title   = {Extensive Chaos in the {L}orenz-96 Model},
  journal = {Chaos},
  year    = {2010},
  volume  = {20},
  number  = {4},
  eid   = {043105},
  nodoi     = {10.1063/1.3496397}
}

@article{katzfuss2016understanding,
  title={Understanding the ensemble {K}alman filter},
  author={Katzfuss, Matthias and Stroud, Jonathan R and Wikle, Christopher K},
  journal={Amer. Statist.},
  volume={70},
  number={4},
  pages={350--357},
  year={2016},
  publisher={Taylor \& Francis}
}

@incollection{lorenz1996predictability,
  author    = {Lorenz, Edward N.},
  title     = {Predictability: A Problem Partly Solved},
  booktitle = {Seminar on Predictability, Vol.\ I},
  editor    = {},
  publisher = {ECMWF, Reading, Berkshire, UK},
  year      = {1996},
  pages     = {1--18}
}

@book{moral2004feynman,
    AUTHOR = {Del Moral, Pierre},
     TITLE = {Feynman-{K}ac formulae: Genealogical and interacting particle systems with
              applications},
    noSERIES = {Probability and its Applications (New York)},
      noNOTE = {Genealogical and interacting particle systems with
              applications},
 PUBLISHER = {Springer-Verlag, New York},
      YEAR = {2004},
     PAGES = {xviii+555},
      noISBN = {0-387-20268-4},
   MRCLASS = {60-02 (47D08 60G35 60J10 60K35 82-02 82C22 82C80 92D25)},
  MRNUMBER = {2044973},
MRREVIEWER = {Frederi\ G.\ Viens},
       noDOI = {10.1007/978-1-4684-9393-1},
       noURL = {https://doi.org/10.1007/978-1-4684-9393-1},
}

@inproceedings{schoenholz2017deep,
  title={Deep Information Propagation},
  author={Schoenholz, Samuel S and Gilmer, Justin and Ganguli, Surya and Sohl-Dickstein, Jascha},
  booktitle={Proc. Int. Conf. Learn. Represent.},
  year={2017}
}

@book {Scott1992,
    AUTHOR = {Scott, David W.},
     TITLE = {Multivariate density estimation},
    noSERIES = {Wiley Series in Probability and Mathematical Statistics:
              Applied Probability and Statistics},
      noNOTE = {Theory, practice, and visualization,
              A Wiley-Interscience Publication},
 PUBLISHER = {John Wiley \& Sons, Inc., New York},
      YEAR = {1992},
     PAGES = {xiv+317},
      noISBN = {0-471-54770-0},
   MRCLASS = {62G05 (62-09 62H12)},
  MRNUMBER = {1191168},
MRREVIEWER = {Theo\ Gasser},
       noDOI = {10.1002/9780470316849},
       noURL = {https://doi.org/10.1002/9780470316849},
}

@ARTICLE{2020SciPy-NMeth-short,
  author  = {{P. Virtanen, {et~al.}}},
  title   = {{SciPy} 1.0: {F}undamental Algorithms for Scientific Computing in {P}ython},
  journal = {Nature Methods},
  year    = {2020},
  volume  = {17},
  number  = {3},
  pages   = {261--272},
  nodoi     = {10.1038/s41592-019-0686-2}
}

@book{silverman1986density,
  title={Density Estimation for Statistics and Data Analysis},
  author={Silverman, BW},
  year={1986},
  publisher={Chapman \& Hall/CRC}
}

@inproceedings{snyder2011particle,
  title={Particle filters, the “optimal” proposal and high-dimensional systems},
  author={Snyder, Chris},
  booktitle={Proceedings of the ECMWF Seminar on Data Assimilation for atmosphere and ocean},
  pages={1--10},
  year={2011}
}

@article{snyder2015performance,
  title={Performance bounds for particle filters using the optimal proposal},
  author={Snyder, Chris and Bengtsson, Thomas and Morzfeld, Mathias},
  journal={Mon. Weather Rev.},
  volume={143},
  nonumber={11},
  pages={4750--4761},
  year={2015}
}

@article{schlogl1972chemical,
  title={Chemical reaction models for non-equilibrium phase transitions},
  author={Schl{\"o}gl, Friedrich},
  journal={Zeitschrift f{\"u}r physik},
  volume={253},
  number={2},
  pages={147--161},
  year={1972},
  publisher={Springer}
}

@article {van2019symmetries,
    AUTHOR = {van Kekem, Dirk L. and Sterk, Alef E.},
     TITLE = {Symmetries in the {L}orenz-96 model},
   JOURNAL = {Internat. J. Bifur. Chaos Appl. Sci. Engrg.},
  FJOURNAL = {International Journal of Bifurcation and Chaos in Applied
              Sciences and Engineering},
    VOLUME = {29},
      YEAR = {2019},
    NUMBER = {1},
     PAGES = {1950008, 18},
      noISSN = {0218-1274,1793-6551},
   MRCLASS = {37C80 (34A33 34C45)},
  MRNUMBER = {3911324},
       noDOI = {10.1142/S0218127419500081},
       noURL = {https://doi.org/10.1142/S0218127419500081},
}

@article{vellela2009stochastic,
  title={Stochastic dynamics and non-equilibrium thermodynamics of a bistable chemical system: the {S}chl{\"o}gl model revisited},
  author={Vellela, Melissa and Qian, Hong},
  journal={J. R. Soc. Interface},
  volume={6},
  number={39},
  pages={925--940},
  year={2009},
  publisher={The Royal Society}
}

@article{Vlysidis2018,
  author = {Vlysidis, Michail and Kaznessis, Yiannis N},
  title = {Solving Stochastic Reaction Networks with Maximum Entropy {L}agrange Multipliers},
  journal = {Entropy},
  volume = {20},
  number = {9},
  pages = {678},
  year = {2018}
}

@article{wilks2005effects,
  author  = {Wilks, Daniel S.},
  title   = {Effects of Stochastic Parametrizations in the {L}orenz-96 System},
  journal = {Q. J. R. Meteorol. Soc.},
  year    = {2005},
  volume  = {131},
  pages   = {389--407},
  nodoi     = {10.1256/qj.04.03}
}

@book{Oksendal,
  title={Stochastic Differential Equations: An Introduction with Applications},
  author={{\O}ksendal, Bernt},
  year={2003},
  publisher={Springer Science \& Business Media}
}
\end{refsection}

\printaddresseshere
%\clearpage

%\bibliographystyle{siamplain}
%\bibliography{references}

\end{document}